\documentclass[10pt,a4paper]{amsart}

\usepackage[T1]{fontenc}
\usepackage{lmodern}
\usepackage{amssymb,mathrsfs}
\usepackage{tikz-cd}
\usepackage[margin=1.1in]{geometry}
\usepackage{microtype}
\usepackage{enumitem}
\usepackage[hidelinks]{hyperref}

\makeatletter
\newlength{\biblabeltextwidth}
\renewcommand{\@biblabel}[1]{\makebox[\biblabeltextwidth][l]{[#1]}}
\makeatother

\setlist[enumerate]{label=\textup{(\roman*)},leftmargin=2.5em,itemsep=2pt,topsep=3pt}

\newtheorem{theorem}{Theorem}[section]
\newtheorem{proposition}[theorem]{Proposition}
\newtheorem{lemma}[theorem]{Lemma}
\newtheorem{corollary}[theorem]{Corollary}
\newtheorem{introtheorem}{Theorem}

\newtheorem{introproposition}[introtheorem]{Proposition}

\theoremstyle{definition}
\newtheorem{definition}[theorem]{Definition}
\newtheorem{convention}[theorem]{Convention}

\theoremstyle{remark}
\newtheorem{remark}[theorem]{Remark}
\newtheorem{warning}[theorem]{Warning}

\numberwithin{equation}{section}

\newcommand{\bbA}{\mathbb A}
\newcommand{\Gm}{\mathbb G_m}
\newcommand{\bbP}{\mathbb P}
\newcommand{\bbQ}{\mathbb Q}
\newcommand{\bbR}{\mathbb R}
\newcommand{\bbZ}{\mathbb Z}

\newcommand{\cD}{\mathcal D}
\newcommand{\cE}{\mathcal E}
\newcommand{\cF}{\mathcal F}
\newcommand{\cH}{\mathcal H}
\newcommand{\cL}{\mathcal L}
\newcommand{\cO}{\mathcal O}
\newcommand{\cQ}{\mathcal Q}
\newcommand{\cV}{\mathcal V}

\newcommand{\an}{\mathrm{an}}
\newcommand{\ad}{\mathrm{ad}}
\newcommand{\alg}{\mathrm{alg}}
\newcommand{\adm}{\mathrm a}
\newcommand{\wa}{\mathrm{wa}}
\newcommand{\ssF}{{\mathrm{ss},F}}
\newcommand{\BdR}{B_{\mathrm{dR}}}

\DeclareMathOperator{\GL}{GL}
\DeclareMathOperator{\Res}{Res}
\DeclareMathOperator{\Gr}{Gr}
\DeclareMathOperator{\Spec}{Spec}
\DeclareMathOperator{\Pic}{Pic}
\DeclareMathOperator{\BB}{BB}
\DeclareMathOperator{\Hom}{Hom}
\DeclareMathOperator{\End}{End}
\DeclareMathOperator{\Rep}{Rep}
\DeclareMathOperator{\Fil}{Fil}
\DeclareMathOperator{\gr}{gr}
\DeclareMathOperator{\supp}{supp}
\DeclareMathOperator{\wt}{wt}
\DeclareMathOperator{\conv}{conv}
\DeclareMathOperator{\id}{id}
\DeclareMathOperator{\pr}{pr}
\DeclareMathOperator{\sep}{sep}

\title{Cocompactness in $p$-adic period domains}
\author{Vladyslav Shram}
\address{Mathematisches Institut, Universit\"at M\"unster, Einsteinstrasse 62, M\"unster, Germany}
\email{vshram@uni-muenster.de}

\begin{document}

\begin{abstract}
Let $F$ be a non-archimedean local field and $L/F$ a complete non-archimedean extension. We show that for any connected reductive $F$-group $H$ and any projective $L$-variety $X$ endowed with an $H_L$-action and an ample rational $H_L$-linearization, the semistable locus $X^\ssF\subset X^\an$ defined by the Hilbert--Mumford inequalities for the $F$-rational one-parameter subgroups of $H$ is cocompact for the action of $H(F)$. As an application, the weakly admissible locus $\cF(G,\mu,b)^\wa$ attached to an arbitrary period datum $(G,\{\mu\},[b])$ over a $p$-adic field $F$ is covered by the $J_b(F)$-translates of one compact subset. Furthermore, for a minuscule period datum, we prove that the admissible locus $\cF(G,\mu,b)^\adm$ is covered by the $J_b(F)$-translates of one compact subset if and only if $\cF(G,\mu,b)^\adm=\cF(G,\mu,b)^\wa$. In particular, this proves \cite[Conjecture~7.2]{CFS} of Chen--Fargues--Shen.
\end{abstract}

\maketitle

\section*{Introduction}

Let $F/\bbQ_p$ be a finite extension, let $G$ be a connected reductive group over $F$, and let $\{\mu\}$ be a $G(\overline F)$-conjugacy class of minuscule cocharacters of $G_{\overline F}$. Fix $b\in G(\breve F)$, where $\breve F$ is the completion of the maximal unramified extension of $F$. Associated with these data are two open analytic subspaces of the flag variety $\cF(G,\mu)$ over $\breve E$,
\[
\cF(G,\mu,b)^\adm\subset\cF(G,\mu,b)^\wa\subset\cF(G,\mu)^\an,
\]
the \emph{admissible} resp. \emph{weakly admissible} locus. Here $\breve E$ is the completion of the maximal unramified extension of the reflex field of $\{\mu\}$, and $\cF(G,\mu)^\an$ is the Berkovich analytification of $\cF(G,\mu)$. The weakly admissible locus $\cF(G,\mu,b)^\wa$ is a vast generalization of the Drinfeld upper half-space and was first defined by Rapoport and Zink \cite{RZ} using a certain semistability condition on filtered isocrystals. This locus is of algebraic nature in the sense that weak admissibility of an analytic point depends only on its support in $\cF(G,\mu)$; moreover, it admits an explicit description as the complement of a profinite union of Zariski-closed Schubert varieties in $\cF(G,\mu)^\an$. By contrast, the existence of the admissible locus generally remained conjectural for many years. Rapoport and Zink predicted the existence of an open subspace $\cF(G,\mu,b)^\adm\subset\cF(G,\mu,b)^\wa$ having the same classical points as $\cF(G,\mu,b)^\wa$ and carrying a $p$-adic local system interpolating the crystalline Galois representations attached to those points \cite[Conjecture~11.4.4]{DOR}. Following earlier constructions in certain PEL cases by Hartl \cite{Hartl08,Hartl13} and Faltings \cite{Faltings}, a general construction became available relatively recently through the theory of $G$-bundles on the Fargues--Fontaine curve, based on work of Fargues--Fontaine \cite{FF}, Fargues \cite{Fargues}, Kedlaya--Liu \cite{KL}, and Scholze--Weinstein \cite{SW}.

Although $\cF(G,\mu,b)^\adm$ and $\cF(G,\mu,b)^\wa$ have the same classical points, they need not coincide as analytic spaces. For basic $[b]\in B(G,\mu)$, Chen--Fargues--Shen proved that the two loci coincide if and only if the pair $(G,\{\mu\})$ is fully Hodge--Newton decomposable \cite[Theorem~6.1]{CFS}. The fully Hodge--Newton decomposable pairs were classified by G\"ortz--He--Nie \cite{GHN}, and the question of when $\cF(G,\mu,b)^\adm=\cF(G,\mu,b)^\wa$ for non-basic $[b]$ was subsequently settled by Chen \cite{Chen}.

Motivated by the analogy with Siegel domains for Hermitian symmetric spaces, Chen--Fargues--Shen conjectured that, for minuscule $\mu$ and basic $[b]\in B(G,\mu)$, the equality $\cF(G,\mu,b)^\adm=\cF(G,\mu,b)^\wa$ is equivalent, in the Berkovich formulation, to the existence of a compact subset $D\subset\cF(G,\mu,b)^\adm$ such that $J_b(F)\cdot D=\cF(G,\mu,b)^\adm$ \cite[Conjecture~7.2]{CFS}. Here $J_b(F)=\{g\in G(\breve F): gb=b\sigma(g)\}$ is the group of automorphisms of the isocrystal with $G$-structure determined by $b$, where $\sigma$ denotes the Frobenius of $\breve F/F$.

The forward implication of the conjecture was previously known in a few special cases. The Lubin--Tate and Drinfeld cases provide two classical examples. Beyond these, it was established by Shen in the unramified unitary case of signature $(1,n-1)$ using the Harder--Narasimhan formalism for finite flat group schemes and the Bruhat--Tits stratification of the corresponding unitary Rapoport--Zink space \cite{Shen}. To the best of our knowledge, neither implication of \cite[Conjecture~7.2]{CFS} was previously known in general.

In this paper we prove this conjecture for arbitrary minuscule period data, without assuming that $[b]$ is basic or even that $[b]\in B(G,\mu)$. In fact, we establish the following two stronger statements.

\begin{introtheorem}[Theorems~\ref{thm:wa-cocompactness} and~\ref{thm:admissible-orbit-closure}]
\label{thm:intro-period-domains}
Let $(G,\{\mu\},[b])$ be an arbitrary period datum over $F$. Then:
\begin{enumerate}
\item There exists a compact subset $D\subset\cF(G,\mu,b)^\wa$ such that
\[
J_b(F)\cdot D=\cF(G,\mu,b)^\wa.
\]
\item If $\mu$ is minuscule, then for every compact subset $C\subset\cF(G,\mu,b)^\adm$, one has
\[
\overline{J_b(F)\cdot C}\cap\cF(G,\mu,b)^\wa \subset\cF(G,\mu,b)^\adm.
\]
\end{enumerate}
\end{introtheorem}

The forward implication of the Chen--Fargues--Shen conjecture is then a special case of assertion (i). Assertion (ii) supplies the converse. Indeed, suppose that a compact subset $D\subset\cF(G,\mu,b)^\adm$ satisfies $J_b(F)\cdot D=\cF(G,\mu,b)^\adm$. Then, taking $C=D$ in (ii), we get
\[
\overline{\cF(G,\mu,b)^\adm}\cap\cF(G,\mu,b)^\wa \subset\cF(G,\mu,b)^\adm.
\]
Since the classical points of $\cF(G,\mu,b)^\wa$ are dense and all belong to $\cF(G,\mu,b)^\adm$, the left-hand side equals $\cF(G,\mu,b)^\wa$. Hence $\cF(G,\mu,b)^\adm=\cF(G,\mu,b)^\wa$.

To prove assertion (i), we establish a general cocompactness theorem for rational Hilbert--Mumford semistable loci in projective varieties over non-archimedean fields. Our general result is formulated in the following setting. Let $F$ be a non-archimedean local field, let $L/F$ be a complete non-archimedean extension, let $H$ be a connected reductive group over $F$, and let $X$ be a projective variety over $L$ endowed with an action of $H_L$. Let $\cL\in\Pic^{H_L}(X)\otimes_\bbZ\bbQ$ be an ample rational $H_L$-linearization on $X$. We denote by $X^\ssF\subset X^\an$ the locus of points $x$ such that
\[
\mu^\cL(\supp_X(x),\lambda)\ge 0
\]
for every one-parameter subgroup $\lambda:\Gm\to H$ defined over $F$. We prove the following theorem.

\begin{introtheorem}[Theorem~\ref{thm:semistable-cocompactness}]
\label{thm:intro-semistable-cocompactness}
There exists a compact subset $D\subset X^\ssF$ such that
\[
H(F)\cdot D = X^\ssF.
\]
\end{introtheorem}

Equivalently, the topological quotient $H(F)\backslash X^\ssF$ is quasicompact. Assertion (i) of Theorem~\hyperref[thm:intro-period-domains]{A} now follows as a special case, using a result of Dat--Orlik--Rapoport \cite[Theorem~9.7.3(ii) and \S9.8]{DOR}, which identifies weak admissibility with rational Hilbert--Mumford semistability for a suitable ample rational linearization of the action of $H_0=\Res_{F/\bbQ_p}J_b$ on $\cF(G,\mu)$. Since $H_0(\bbQ_p)=J_b(F)$, applying Theorem~\hyperref[thm:intro-semistable-cocompactness]{B} over $\bbQ_p$ gives (i).

We now briefly explain the idea of the proof of Theorem~\hyperref[thm:intro-semistable-cocompactness]{B}. Replacing $\cL$ by a sufficiently large positive power gives an $H_L$-equivariant embedding $X\hookrightarrow\bbP(V)$, where $V=H^0(X,\cL)^\vee$ carries the algebraic $H_L$-action coming from the equivariant structure on $\cL$. Fix a non-archimedean $L$-norm $\|\cdot\|$ on $V$ that is split in some basis. For $x\in X^\an$ represented by a nonzero vector $v\in V_{\cH(x)}$, put
\[
\delta_h(x) = \frac{\|hv\|}{\|v\|},\qquad \delta(x)=\inf_{h\in H(F)} \delta_h(x) = \frac{\inf_{h \in H(F)} \|hv\|}{\|v\|}\in[0,1].
\]
Note that $\delta(hx)=\delta(x)\delta_h(x)^{-1}$ for any $h\in H(F)$. We prove the following orbit-closure criterion for $F$-rational semistability (Theorem~\ref{thm:orbit-closure-criterion}):
\[
x\in X^\ssF \quad\Longleftrightarrow\quad \delta(x)>0 \quad\Longleftrightarrow\quad 0\notin\overline{H(F)\cdot v} \subset V_{\cH(x)},
\]
where the closure is taken in the non-archimedean topology of $V_{\cH(x)}$. Moreover, we show that whenever these conditions hold, the infimum defining $\delta(x)$ is attained.

The main step in the proof of the criterion is the implication from $F$-semistability to $\delta(x)>0$, which requires a positive lower bound for $\|hv\|/\|v\|$ that is uniform in $h\in H(F)$. We obtain such a bound using a decomposition $H(F)=C_1S(F)C_2$, where $S$ is a maximal $F$-split torus and $C_1,C_2 \subset H(F)$ are compact subsets. Using the convex-hull criterion for $S$-semistability (Lemma~\ref{lem:convex-hull-criterion}), we construct a continuous function that is positive at every $S$-semistable point and provides a lower bound for $\|av\|/\|v\|$ for all $a\in S(F)$. Compactness of $C_2x$, together with uniform bounds for the two compact factors, then gives the required estimate over all of $H(F)$.

The orbit-closure criterion now implies Theorem~\hyperref[thm:intro-semistable-cocompactness]{B}. Indeed, for $x\in X^\ssF$, choose $h\in H(F)$ for which $\delta(x)=\delta_h(x)$. Then $\delta(hx)=\delta(x)\delta_h(x)^{-1}=1$, and hence the subset
\[
D=\{x\in X^\an:\delta(x)=1\}
\]
has the desired property $H(F)\cdot D= X^\ssF$. Since $\delta$ is the infimum of continuous functions $\delta_h$ and $\delta\le 1$, the set $D$ is closed in the compact space $X^\an$, and hence compact.

Let us now turn to assertion (ii) of Theorem~\hyperref[thm:intro-period-domains]{A}. Let $C\subset\cF(G,\mu,b)^\adm$ be compact, and suppose that $h_nx_n\to x'\in\cF(G,\mu,b)^\wa$, where $x_n\in C$ and $h_n\in J_b(F)$. Since $\cF(G,\mu)^\an$ is metrizable \cite[Remark~1.5]{HLP}, compactness of $C$ allows us, after passing to a subsequence, to assume that $x_n\to x\in C\subset\cF(G,\mu,b)^\adm$. We then show that weak admissibility of both limits $x$ and $x'$ allows us to relate them by opposite one-parameter subgroups in the following sense. As with Theorem~\hyperref[thm:intro-semistable-cocompactness]{B}, we establish this in the setting of rational Hilbert--Mumford semistability in arbitrary projective varieties.

\begin{introtheorem}[Theorem~\ref{thm:common-specializations}]
\label{thm:intro-common-specializations}
In the setting of Theorem~\hyperref[thm:intro-semistable-cocompactness]{B}, suppose that
\[
x_n\longrightarrow x,\qquad h_nx_n\longrightarrow x',
\]
where $x_n\in X^\an$, $h_n\in H(F)$ and $x,x'\in X^\ssF$. Then there exist an $F$-rational one-parameter subgroup $\lambda:\Gm\to H$ and an element $a\in H(F)$ such that
\[
\lim_{t\to0}\lambda(t)x = \lim_{t\to0}\lambda(t)^{-1}ax' \in X^\ssF.
\]
\end{introtheorem}

The key idea in the proof of Theorem~\hyperref[thm:intro-common-specializations]{C} is that, after reducing to a maximal $F$-split torus $S$ by the Cartan decomposition, the two convergences $x_n\to x$ and $h_nx_n\to x'$ impose strong restrictions on the $S$-weight components of $x$ and $x'$. Together with semistability, these restrictions allow us to find an $F$-rational cocharacter $\lambda$ which acts with positive weights on the $S$-weight components occurring only in $x$, with negative weights on those occurring only in $x'$, and with weight zero on the common components. Taking the limits $\lim_{t\to0}\lambda(t)x$ and $\lim_{t\to0}\lambda(t)^{-1}x'$ then amounts to projecting $x$ and $x'$ onto the $\lambda$-weight-zero subspace of $V$, and comparing these projections shows that the two limits agree up to translation by an element of $S(F)$.

Applied to our period-domain setting, Theorem~\hyperref[thm:intro-common-specializations]{C} gives an $F$-rational cocharacter $\lambda:\Gm\to J_b$ and an element $a\in J_b(F)$ such that
\[
y:=\lim_{t\to0}\lambda(t)x= \lim_{t\to0}\lambda(t)^{-1}ax'\in\cF(G,\mu,b)^\wa.
\]
Now, in order to prove that $x'$ is admissible, it is enough to show that $y$ is admissible. Indeed, if $y\in\cF(G,\mu,b)^\adm$, then $\lambda(t)^{-1}ax'\in\cF(G,\mu,b)^\adm$ for all small enough $t\in F^\times$. Since the admissible locus is $J_b(F)$-invariant, this gives $x'\in\cF(G,\mu,b)^\adm$. Thus the proof of assertion (ii) of Theorem~\hyperref[thm:intro-period-domains]{A} has now been reduced to the following statement.

\begin{introproposition}[Proposition~\ref{prop:admissible-1PS-specialization}]
\label{prop:intro-admissible-specialization}
Assume that $\mu$ is minuscule. Let $x\in \cF(G,\mu,b)^\adm$ and let $\lambda:\Gm\to J_b$ be a one-parameter subgroup over $F$. Let $x_0 = \lim_{t\to 0}\lambda(t)x\in \cF(G,\mu)^\an$. If $x_0\in \cF(G,\mu,b)^\wa$, then $x_0 \in \cF(G,\mu,b)^\adm$.
\end{introproposition}

Proposition~\hyperref[prop:intro-admissible-specialization]{D} is the only ingredient in the proof of Theorem~\hyperref[thm:intro-period-domains]{A} that is specific to $p$-adic period domains and whose proof uses modifications of bundles on the Fargues--Fontaine curve. After passing to a representation of $G$, we can use the description of the limit $\lim_{t\to 0}\lambda(t)x$ in the flag variety as the associated graded filtration to identify the modified vector bundle corresponding to $x_0$ with the direct sum of subquotients of the modified bundle corresponding to $x$. The weak admissibility of $x_0$ then forces all those subquotients to have degree $0$. Finally, admissibility of $x$ allows us to deduce that each of those subquotients, and hence the modified bundle corresponding to $x_0$, is trivial.

The paper is organized as follows. In Section~\ref{sec:rational-semistability} we introduce rational Hilbert--Mumford semistability, prove the orbit-closure criterion, and deduce Theorem~\hyperref[thm:intro-semistable-cocompactness]{B} from it. In Section~\ref{sec:common-specializations} we prove Theorem~\hyperref[thm:intro-common-specializations]{C}. Section~\ref{sec:weakly-admissible-period-domains} then exhibits the weakly admissible locus as an instance of the rational Hilbert--Mumford semistable locus and deduces Theorem~\hyperref[thm:intro-period-domains]{A}(i). Finally, in Section~\ref{sec:admissible-period-domains} we prove Proposition~\hyperref[prop:intro-admissible-specialization]{D} and deduce Theorem~\hyperref[thm:intro-period-domains]{A}(ii).

\subsection*{Acknowledgements}
I am deeply grateful to my advisor Eva Viehmann for suggesting that I think about the cocompactness conjecture, for her guidance and support at various stages of this project, and for the many mathematical discussions that have contributed to this work. Furthermore, I would like to thank Eugen Hellmann for his questions and comments, and for pointing out a mistake in an earlier approach to the problem. I also thank Zhixiang Wu and Ioannis Zachos for helpful discussions.

The author was supported by the DFG through the CRC 1442 ``Geometry: Deformations and Rigidity'' (project number 427320536), under Germany's Excellence Strategy EXC 2044/2--390685587 ``Mathematics M\"unster: Dynamics--Geometry--Structure'', and through the CRC/TRR 326 ``Geometry and Arithmetic of Uniformized Structures'' (project number 444845124), as well as by the Leibniz Prize of Eva Viehmann.

\section{Cocompactness of $F$-semistable loci}\label{sec:rational-semistability}
\subsection{$F$-rational Hilbert--Mumford semistability}

Throughout Sections~\ref{sec:rational-semistability}--\ref{sec:common-specializations}, let $F$ be a non-archimedean local field and $L/F$ a complete non-archimedean extension. Let $H$ be a connected reductive group over $F$. Let $X$ be a projective variety over $L$ equipped with an action of $H_L$, and let $\cL\in \Pic^{H_L}(X)\otimes_\bbZ\bbQ$ be an ample rational $H_L$-equivariant line bundle; we shall also refer to $\cL$ as a \emph{rational $H_L$-linearization}. Let $X^\an$ denote the Berkovich analytification of $X$ over the valued field $L$; we use the same notation for its underlying topological space. We have a natural surjective continuous map $\supp_X:X^\an\to X$, which affine-locally is given by sending a multiplicative seminorm to its support prime ideal. Since $X$ is projective, $X^\an$ is compact Hausdorff \cite[Theorem~3.4.8(ii)]{Berkovich}.

We recall the Hilbert--Mumford numerical function. Let $\xi\in X$, viewed as its canonical $\kappa(\xi)$-valued point, and let $\lambda:\Gm\to H$ be a one-parameter subgroup defined over $F$. Since $X$ is proper, the orbit map
\[
\mathbb G_{m,\kappa(\xi)}\longrightarrow X_{\kappa(\xi)}, \quad t\longmapsto\lambda(t)\xi
\]
extends uniquely to $\bbA^1_{\kappa(\xi)}$. Write $\xi_0=\lim_{t\to0}\lambda(t)\xi$ for its value at $0$. Let $\mathcal M\in\Pic^{H_L}(X)$ be an $H_L$-equivariant line bundle on $X$, i.e.~a line bundle equipped with an isomorphism $\phi:\alpha^*\mathcal M\xrightarrow{\sim}p_2^*\mathcal M$ satisfying the usual cocycle condition on $H_L\times H_L\times X$, where $\alpha:H_L\times X\to X$ is the action map and $p_2:H_L\times X\to X$ is the second projection. Since $\xi_0$ is fixed by $\lambda$, evaluating $\phi$ at $(\lambda(t),\xi_0)$ for $t\in\mathbb G_{m,\kappa(\xi)}$ defines an algebraic action of $\mathbb G_{m,\kappa(\xi)}$ on the fiber $\mathcal M(\xi_0)$. This action is given by a character $t\mapsto t^r$ for some $r\in\bbZ$, and we put $\mu^{\mathcal M}(\xi,\lambda):=-r$.

For $\cL\in \Pic^{H_L}(X)\otimes_\bbZ\bbQ$, choose a positive integer $N$ such that $N\cL \in \Pic^{H_L}(X)$, and set $\mu^\cL(\xi,\lambda) = \frac{1}{N}\mu^{N\cL}(\xi,\lambda)$. The additivity of Hilbert--Mumford weights under tensor products shows that this is independent of the choice of $N$; see \cite[\S2.2]{DOR}.

\begin{definition}\label{def:rational-semistability}
A point $\xi\in X$ is called \emph{$F$-rationally semistable}, or simply \emph{$F$-semistable}, if, viewed as the canonical $\kappa(\xi)$-valued point of $X$, it satisfies $\mu^\cL(\xi,\lambda)\ge 0$ for every one-parameter subgroup $\lambda:\Gm\to H$ defined over $F$. Let $X^{\ssF,\alg}\subset X$ denote the locus of $F$-semistable points. We then define the $F$-semistable locus in $X^\an$ by
\[
X^\ssF:=\supp_X^{-1}(X^{\ssF,\alg})\subset X^\an.
\]
\end{definition}

In particular, the semistability of a point of $X^\an$ depends only on its support; for $x\in X^\an$, we also write $\mu^\cL(x,\lambda):=\mu^\cL(\supp_X(x),\lambda)$. Moreover, $X^{\ssF,\alg}$, and hence also $X^\ssF$, is stable under the action of $H(F)$, since $\mu^\cL(h\xi,\lambda)=\mu^\cL(\xi,h^{-1}\lambda h)$ and $h^{-1}\lambda h$ is defined over $F$ whenever $\lambda$ is.

The notion of semistability defined above depends on the fields $F$ and $L$, the group $H$, its action on $X$, and the linearization $\cL$. We usually leave this dependence implicit in our terminology and notation when these data are clear from the context. When we wish to emphasize the group $H$ with respect to which semistability is considered, we may explicitly say that a point is \emph{$H$-semistable}.

Since Hilbert--Mumford weights are additive under tensor products, replacing $\cL$ by a positive multiple does not change the $F$-semistable locus. We may therefore assume that $\cL$ is a very ample $H_L$-equivariant line bundle. Let
\[
\iota: X\hookrightarrow\bbP(V)
\]
be the corresponding $H_L$-equivariant embedding, where $V=H^0(X,\cL)^\vee$ carries an algebraic action $\rho:H_L\to\GL(V)$ coming from the $H_L$-equivariant structure on $\cL$. We have $\iota^*\cO(1)\cong\cL$ as $H_L$-equivariant line bundles.

\begin{remark}\label{rem:zariski-non-openness}
In ordinary geometric invariant theory, the semistable locus is Zariski open. In contrast, the locus $X^{\ssF,\alg}$ obtained by testing only the one-parameter subgroups defined over $F$ need not be Zariski open. For example, in the Drinfeld case with $n\ge2$, it is the complement in $\bbP_F^{n-1}$ of the union of all $F$-rational hyperplanes. In particular, it contains the generic point but no $F$-rational point and is therefore not Zariski open, since $\bbP_F^{n-1}(F)$ is Zariski dense in $\bbP_F^{n-1}$. Nevertheless, the inverse image $X^\ssF = \supp_X^{-1}(X^{\ssF,\alg})$ is open in $X^\an$, as will be shown in Corollary~\ref{cor:semistability-openness}.
\end{remark}

\begin{convention}\label{conv:projective-representative}
Let $x\in X^\an$ be a point, let $\xi=\supp_X(x)\in X$ be its support, and let $K=\cH(x)$ be the completed residue field of $x$. Then $K/L$ is a complete valued field extension, and $x$ determines a $K$-valued point $x_K:\Spec K\to X$ whose image is $\xi$. Put $V_K=V\otimes_LK$. After the fixed embedding $\iota: X\hookrightarrow\bbP(V)$, choose a nonzero vector $v\in V_K$ representing this point, so that $x_K=[v]\in X(K)\subset\bbP(V)(K)$. We shall then say that $x$ is \textit{represented} by $v \in V_K$.
\end{convention}

Let $x\in X^\an$ be represented by $v\in V_K$ and let $\lambda:\Gm\to H$ be a one-parameter subgroup defined over $F$. Let $V=\bigoplus_{m\in\bbZ}V_m$ be the weight decomposition with respect to the action of $\lambda$. Then, writing $v=\sum_m v_m$ with $v_m\in V_{m,K}$, we have
\begin{equation}\label{eq:hilbert-mumford-weight-projective-formula}
\mu^\cL(x,\lambda) = \mu^{\cO(1)}([v],\lambda) = -\min\{m:v_m\ne 0\};
\end{equation}
see e.g.~\cite[\S2.2]{DOR}.

Let $S\subset H$ be a maximal $F$-split torus, and let $V=\bigoplus_{\chi\in X^*(S)}V_\chi$ be the corresponding relative weight decomposition. For a point $x \in X^\an$ represented by $v \in V_K$, write $v=\sum_\chi v_\chi$ according to this weight decomposition, and put
\[
\wt_S(x)=\{\chi\in X^*(S):v_\chi\ne0\}.
\]
We regard this finite set as a subset of $X^*(S)_\bbR$ and write $\conv(\wt_S(x))$ for its convex hull. Using formula \eqref{eq:hilbert-mumford-weight-projective-formula}, one obtains the following standard weight-polytope criterion of semistability.

\begin{lemma}\label{lem:convex-hull-criterion}
Let $x\in X^\an$, and let $S\subset H$ be a maximal $F$-split torus. Then $x$ is $S$-semistable if and only if
\[
0\in\conv(\wt_S(x))\subset X^*(S)_\bbR.
\]
\end{lemma}
\begin{proof}
See e.g.~\cite[Theorem~9.2]{Dolgachev}.
\end{proof}

Note that a point $x \in X^\an$ is $H$-semistable if and only if it is $S$-semistable for all maximal $F$-split tori $S \subset H$.

\subsection{The orbit-closure criterion for $F$-semistability}

For the rest of this section and Section~\ref{sec:common-specializations}, fix a maximal $F$-split torus $S \subset H$. Let $\wt_S(V) = \{\chi \in X^*(S): V_\chi\neq0\}$ be the set of $S$-weights of $V$, and let $V = \bigoplus_{\chi \in \wt_S(V)}V_\chi$ be the corresponding weight decomposition. Choose an $L$-basis $(e_{\chi,j})_j$ of each $V_\chi$, so that their union forms an $L$-basis of $V$, and fix positive real numbers $b_{\chi,j}$. For every complete non-archimedean extension $K/L$, equip $V_K$ with the non-archimedean $K$-norm
\[
\Bigl\|\sum_{\chi,j}a_{\chi,j}e_{\chi,j}\Bigr\|_K:= \max_{\chi,j}b_{\chi,j}|a_{\chi,j}|,
\]
where $a_{\chi,j}\in K$; in the following, we will simply write $\|\cdot\|$ instead of $\|\cdot\|_K$, suppressing $K$ from the notation. Since the basis is adapted to the $S$-weight decomposition, this norm is split by that decomposition, i.e.
\[
\Bigl\|\sum_\chi v_\chi\Bigr\|=\max_\chi\|v_\chi\|
\]
for $v_\chi \in V_{\chi,K}$.

The algebraic action of $H_L$ on $V$ induces a continuous action of $H(F)$ on $V_K$ for the topologies coming from $F$ and $K$, respectively. For $x\in X^\an$ represented by $v\in V_K$, put
\[
\delta_h(x):= \frac{\|hv\|}{\|v\|},\qquad \delta(x):=\inf_{h\in H(F)} \delta_h(x) = \frac{\inf_{h\in H(F)} \|hv\|}{\|v\|}\in[0,1].
\]
Thus $\delta(hx) = \delta(x)\delta_h(x)^{-1}$ for any $h \in H(F)$ and $x \in X^\an$. Although each $\delta_h$ is continuous on $X^\an$, their infimum $\delta$ need not be continuous. For each $\chi\in\wt_S(V)$, put
\begin{equation}\label{eq:chi-ratio}
r_\chi(x):=\frac{\|v_\chi\|}{\|v\|}.
\end{equation}
Each $r_\chi$ is a well-defined continuous function on $X^\an$, and $\wt_S(x)=\{\chi\in\wt_S(V):r_\chi(x)>0\}$. Let $\mathscr Q_S=\{\varnothing\ne Q\subset\wt_S(V):0\in\conv(Q)\}$, and define
\[
\varepsilon_S(x) := \max_{Q\in\mathscr Q_S} \min_{\chi\in Q} r_\chi(x),
\]
with the convention that $\varepsilon_S=0$ when $\mathscr Q_S=\varnothing$. As the maximum and minimum are taken over finite sets, $\varepsilon_S$ is continuous on $X^\an$. Moreover, $\varepsilon_S(x)>0$ if and only if $0\in\conv(\wt_S(x))$, which, by Lemma~\ref{lem:convex-hull-criterion}, is precisely the condition that $x$ be $S$-semistable.

By the Cartan decomposition \cite[\S1.2]{Cartier}, there are compact subsets $C_1,C_2\subset H(F)$ such that
\begin{equation}\label{eq:cartan-decomposition}
H(F)=C_1S(F)C_2.
\end{equation}
In view of this decomposition, put
\[
\varepsilon_S^{C_2}(x):=\min_{c\in C_2}\varepsilon_S(cx).
\]
The continuity of the action $H(F) \times X^\an \to X^\an$ and the compactness of $C_2$ show that $\varepsilon_S^{C_2}$ is a well-defined continuous function on $X^\an$. The following estimate relates it to $\delta$.

\begin{lemma}\label{lem:delta-lower-bound}
There exists a constant $\gamma>0$ such that
\[
\delta(x)\ge\gamma\varepsilon_S^{C_2}(x)
\]
for every $x\in X^\an$.
\end{lemma}
\begin{proof}
Let $\|\cdot\|_{\mathrm{op}}$ be the operator norm on $\End_L(V)$ induced by our fixed norm $\|\cdot\|$ on $V$. The function $c\mapsto\|\rho(c)^{-1}\|_{\mathrm{op}}$ is continuous on $H(F)$, and therefore attains a maximum on each of the compact sets $C_1$ and $C_2$. Put
\[
\gamma_i:=\left(\max_{c\in C_i}\|\rho(c)^{-1}\|_{\mathrm{op}}\right)^{-1}>0,\qquad i=1,2.
\]
Take any $c\in C_i$ and a complete non-archimedean extension $K/L$. Let $e_1,\dots,e_N$ be our fixed $L$-basis of $V$ from above. For any $v\in V_K$, write $cv=\sum_j a_je_j$, with $a_j\in K$. Since $\|c^{-1}e_j\|\le\gamma_i^{-1}\|e_j\|$ for every $j$, we have
\[
\|v\|=\Bigl\|\sum_j a_jc^{-1}e_j\Bigr\|\le \max_j|a_j|\,\|c^{-1}e_j\|\le \gamma_i^{-1}\max_j|a_j|\,\|e_j\|=\gamma_i^{-1}\|cv\|.
\]
This shows that $\|cv\|\ge\gamma_i\|v\|$ for all $c\in C_i$ and $v\in V_K$, with the same constant $\gamma_i$ for every complete non-archimedean extension $K/L$.

Now let $x\in X^\an$ be represented by $v\in V_K$, with $K=\cH(x)$, and let $a\in S(F)$. For every $Q\in\mathscr Q_S$, the condition $0\in\conv(Q)$ implies that the $\bbR$-linear map
\[
X^*(S)_\bbR\longrightarrow\bbR, \quad \chi\longmapsto\log|\chi(a)|
\]
cannot be strictly negative on all of $Q$. Thus $\max_{\chi\in Q}|\chi(a)|\ge 1$, and hence
\[
\|av\|=\max_{\chi\in\wt_S(V)}|\chi(a)|\,\|v_\chi\|\ge\max_{\chi\in Q}|\chi(a)|\,\|v_\chi\| =\left(\max_{\chi\in Q}|\chi(a)|\,r_\chi(x)\right)\|v\|\ge\left(\min_{\chi\in Q}r_\chi(x)\right)\|v\|.
\]
By the definition of $\varepsilon_S$, we obtain $\|av\|\ge\varepsilon_S(x)\|v\|$ for all $x \in X^\an$ represented by $v \in V_K$.

To pass from $S(F)$ to $H(F)$, let $h\in H(F)$ and write $h=c_1ac_2$ with $c_i\in C_i$ and $a\in S(F)$. For any $x \in X^\an$ represented by $v \in V_K$, the preceding inequalities give
\[
\|hv\|\ge\gamma_1\|ac_2v\|\ge\gamma_1\varepsilon_S(c_2x)\|c_2v\|\ge\gamma_1\gamma_2\varepsilon_S^{C_2}(x)\|v\|.
\]
Taking the infimum over $h\in H(F)$ proves the lemma with $\gamma=\gamma_1\gamma_2$.
\end{proof}

For $x\in X^\an$, the condition $x\in X^\ssF$ is purely algebraic: it depends only on the scheme-theoretic support of $x$ and the $F$-rational one-parameter subgroups of $H$, and is defined without any reference to the non-archimedean topologies of $F$ and $L$, or to the completed residue field $K=\cH(x)$ of $x$. Interestingly, this condition is nevertheless equivalent to one formulated entirely in terms of the non-archimedean topology of $V_K$.

\begin{theorem}\label{thm:orbit-closure-criterion}
Let $x\in X^\an$ be represented by $v\in V_K$. Then the following are equivalent:
\begin{enumerate}
\item $x\in X^\ssF$;
\item $\delta(x)>0$;
\item $0\notin\overline{H(F)\cdot v}$ in $V_K$.
\end{enumerate}
Moreover, when these conditions hold, the infimum defining $\delta(x)$ is attained.
\end{theorem}
\begin{proof}
Conditions (ii) and (iii) are equivalent by the definition of $\delta$. We now prove the equivalence of (i) and (ii).

If $x$ is not semistable, there is an $F$-rational cocharacter $\lambda$ with $\mu^\cL(x,\lambda)<0$. By~\eqref{eq:hilbert-mumford-weight-projective-formula}, this means that $v_m= 0$ for $m\le 0$, so for a uniformizer $\varpi$ of $F$ we have
\[
\lambda(\varpi^n)v\longrightarrow 0, \quad n \to +\infty,
\]
which implies $\inf_{h \in H(F)}\|hv\|=0$, and hence $\delta(x)=0$.

Conversely, assume that $x$ is semistable. Let $S \subset H$ be the maximal $F$-split torus fixed above and let $C_1, C_2 \subset H(F)$ be the compact subsets from \eqref{eq:cartan-decomposition}, so that $H(F)=C_1S(F)C_2$. Since the action of $H(F)$ preserves semistability, every $cx$ with $c \in C_2$ is $H$-semistable and hence also $S$-semistable. In particular, $\varepsilon_S(cx)>0$ for all $c \in C_2$. As $C_2$ is compact and $\varepsilon_S$ is continuous, we obtain $\varepsilon_S^{C_2}(x) = \min_{c \in C_2}\varepsilon_S(cx)>0$. Lemma~\ref{lem:delta-lower-bound} now gives $\delta(x)\ge\gamma \varepsilon_S^{C_2}(x)>0$, proving the equivalence of (i) and (ii).

It remains to prove the attainment assertion. Let $x\in X^\an$ be represented by $v\in V_K$, and assume that it satisfies the equivalent conditions above. The algebraicity of the representation $\rho: H_L \to \GL(V)$ implies that the orbit $H(F)\cdot v$ is contained in a finite-dimensional $F$-subspace of $V_K$. To see this, let $e_1,\dots,e_N$ be the fixed $L$-basis of $V$ and, with respect to this basis, let $x_{ij}$ denote the standard coordinate functions on $\GL(V)$. Put
\[
a_{ij}:=\rho^*(x_{ij})\in \Gamma(H_L,\cO_{H_L})=L\otimes_F\Gamma(H,\cO_H).
\]
Then, writing $a_{ij}=\sum_r c_{ijr}\otimes f_{ijr}$ for $c_{ijr}\in L$, $f_{ijr}\in\Gamma(H,\cO_H)$, and $v=\sum_j v_je_j$ for $v_j \in K$, we have
\[
\rho(h)v = \sum_{i,j,r} f_{ijr}(h)c_{ijr}v_je_i
\]
for any $h \in H(F)$. Since $f_{ijr}(h)\in F$, it follows that
\[
H(F)\cdot v\subset \operatorname{span}_F \{c_{ijr}v_je_i\}_{i,j,r} =: W_v,
\]
which is finite-dimensional over $F$.

Choose $h_n\in H(F)$ such that $\|h_nv\|\to \inf_{h \in H(F)}\|hv\|$. Then $(h_nv)$ is a bounded sequence in $W_v$, which is finite-dimensional over the locally compact field $F$. It therefore has a convergent subsequence, and after passing to it we may assume that $h_nv\to v'\in W_v$. Continuity of the norm gives $\|v'\|=\inf_{h \in H(F)}\|hv\|>0$. Now, since $h_nv \to v'$, we in particular have $\|h_nv-v'\|<\|v'\|$ for all $n\gg 0$, and the ultrametric inequality therefore gives $\|h_nv\|=\|v'\|=\inf_{h \in H(F)}\|hv\|$.
\end{proof}

\begin{corollary}\label{cor:semistability-openness}
The subset $X^\ssF\subset X^\an$ is open.
\end{corollary}
\begin{proof}
As we saw in the proof of Theorem~\ref{thm:orbit-closure-criterion}, every $x\in X^\ssF$ satisfies $\varepsilon_S^{C_2}(x)>0$. Conversely, positivity of $\varepsilon_S^{C_2}(x)$ implies $\delta(x)>0$ by Lemma~\ref{lem:delta-lower-bound}, and hence $x\in X^\ssF$ by Theorem~\ref{thm:orbit-closure-criterion}. It follows that $X^\ssF=\{x\in X^\an:\varepsilon_S^{C_2}(x)>0\}$. Since $\varepsilon_S^{C_2}$ is continuous, this set is open.
\end{proof}

\subsection{Cocompactness of the $F$-semistable locus}

From Theorem~\ref{thm:orbit-closure-criterion} we now deduce a general cocompactness theorem for $F$-semistable loci in projective varieties.

\begin{theorem}\label{thm:semistable-cocompactness}
There exists a compact subset $D\subset X^\ssF$ such that
\[
H(F)\cdot D = X^\ssF.
\]
\end{theorem}
\begin{proof}
Take $D=\{x \in X^\an: \delta(x)=1\}$, which is contained in $X^\ssF$ by Theorem~\ref{thm:orbit-closure-criterion}. For any $x\in X^\ssF$, the attainment assertion of Theorem~\ref{thm:orbit-closure-criterion} gives an $h \in H(F)$ such that $\delta(x)=\delta_h(x)$. Consequently $\delta(hx)=\delta(x)\delta_h(x)^{-1} = 1$, so $hx \in D$. This shows that $H(F)\cdot D = X^\ssF$.

To prove compactness of $D$, we need to show that $D$ is closed in the compact space $X^\an$. Since $\delta=\inf_{h \in H(F)}\delta_h$ and each $\delta_h$ is continuous on $X^\an$, the equality
\[
D = \bigcap_{h \in H(F)}\{x \in X^\an: \delta_h(x)\ge 1\}
\]
exhibits $D$ as an intersection of closed subsets of $X^\an$.
\end{proof}

\begin{remark}[Adic version]\label{rem:adic-semistable-cocompactness}
Let $X^\ad$ denote the adic analytification of $X$ over the valued field $L$. As in the Berkovich case, define the $F$-semistable locus in $X^\ad$ by
\[
X^{\ssF,\ad}:= (\supp_X^\ad)^{-1}(X^{\ssF,\alg}),
\]
where $\supp_X^\ad:X^\ad\to X$ is the analogous natural map of topological spaces, given affine-locally by sending a valuation to its support prime ideal. The adic--Berkovich comparison gives a quotient map of topological spaces $\sep_X:X^\ad\to X^\an$, which is functorial in $X$ and satisfies $\supp_X \circ \sep_X = \supp_X^\ad$ \cite[(8.1.6), Lemma~8.1.7(ii) and Remark~8.3.2]{Huber}. It follows that $X^{\ssF,\ad}=\sep_X^{-1}(X^\ssF)$, and in particular $X^{\ssF,\ad}$ is open in $X^\ad$. By functoriality, $\sep_X$ is also $H(F)$-equivariant, so
\[
H(F)\cdot \sep_X^{-1}(D) = \sep_X^{-1}(H(F)\cdot D) = \sep_X^{-1}(X^\ssF) = X^{\ssF,\ad}.
\]
Since $D$ is closed in $X^\an$, its inverse image $\sep_X^{-1}(D)$ is closed in the quasicompact space $X^\ad$ and is therefore quasicompact. Since $X^{\ssF,\ad}$ is a locally spectral space, $\sep_X^{-1}(D)$ is contained in a quasicompact open subset $U \subset X^{\ssF,\ad}$. The preceding equality then implies $H(F)\cdot U = X^{\ssF,\ad}$, giving the adic version of Theorem~\ref{thm:semistable-cocompactness}.
\end{remark}

\section{Relating semistable orbit limits by one-parameter subgroups}\label{sec:common-specializations}

We retain the general setting of Section~\ref{sec:rational-semistability}. In this section we show that any two semistable limits of sequences related by elements of $H(F)$ admit, after a translation by an element of $H(F)$, common specializations by opposite one-parameter subgroups over $F$. In Section~\ref{sec:admissible-period-domains}, this will be applied to $p$-adic period domains in order to reduce the proof of Theorem~\hyperref[thm:intro-period-domains]{A}(ii) to a question about the behavior of modifications of bundles on the Fargues--Fontaine curve under one-parameter specializations.

For $x\in X^\an$ and an $F$-rational cocharacter $\lambda:\Gm\to H$, put $K=\cH(x)$ and denote by $x_K\in X(K)$ the canonical $K$-valued point associated with $x$. The algebraic one-parameter limit
\[
x_{0,K}:=\lim_{t\to0}\lambda(t)x_K\in X(K)
\]
exists by properness of $X$. We may also regard this as a usual topological limit in $X(K)$, where $X(K) \subset \bbP(V)(K)$ carries the natural topology induced by $K$. Since the map $X(K) \to X^\an$ is continuous, the image $x_0$ of $x_{0,K}$ under this map satisfies $\lambda(t)x\to x_0$ in $X^\an$ as $t\to0$ in $F^\times$, i.e.~$x_0=\lim_{t\to0}\lambda(t)x$.

\begin{theorem}\label{thm:common-specializations}
Suppose that
\[
x_n\longrightarrow x,\qquad h_nx_n\longrightarrow x'
\]
in $X^\an$, where $x_n \in X^\an$, $h_n\in H(F)$ and $x,x'\in X^\ssF$. Then there exist an $F$-rational one-parameter subgroup $\lambda:\Gm\to H$ and $a\in H(F)$ such that
\[
\lim_{t\to0}\lambda(t)x = \lim_{t\to0}\lambda(t)^{-1}ax'
\]
and $\lim_{t\to0}\lambda(t)x \in X^\ssF$.
\end{theorem}

Before proving the theorem, let us record the following consequence of Theorem~\ref{thm:orbit-closure-criterion}, which will be used to deduce semistability of the limit.

\begin{lemma}\label{lem:semistable-1PS-specialization}
Let $x\in X^\ssF$ and let $\lambda:\Gm\to H$ be defined over $F$. Put $x_0=\lim_{t\to0}\lambda(t)x \in X^\an$. Then $x_0 \in X^\ssF$ if and only if $\mu^\cL(x,\lambda)=0$.
\end{lemma}
\begin{proof}
Suppose first that $\mu^\cL(x,\lambda)=0$. Put $K=\cH(x)$, represent $x$ by $v\in V_K$, and write $v=\sum_m v_m$ according to the $\lambda$-weight decomposition. Formula~\eqref{eq:hilbert-mumford-weight-projective-formula} gives $v_m=0$ for $m<0$, and hence, for a uniformizer $\varpi$ of $F$,
\[
\lambda(\varpi^n)v\longrightarrow v_0\ne0
\]
in $V_K$. In particular $v_0\in\overline{H(F)\cdot v}$, and since this set is $H(F)$-stable and closed, we obtain $\overline{H(F)\cdot v_0}\subset \overline{H(F)\cdot v}$. By Theorem~\ref{thm:orbit-closure-criterion}, $0\notin\overline{H(F)\cdot v}$, so also $0\notin\overline{H(F)\cdot v_0}$. But since $v_0\neq0$, it represents $x_0$, and hence another application of Theorem~\ref{thm:orbit-closure-criterion} gives $x_0\in X^\ssF$.

Conversely, suppose that $x_0\in X^\ssF$. Semistability of $x_0$ in particular gives $\mu^\cL(x_0,\lambda)\ge 0$ and $\mu^\cL(x_0,-\lambda)\ge 0$. On the other hand, $x_0$ is fixed by $\lambda$, so $\mu^\cL(x_0,-\lambda)=-\mu^\cL(x_0,\lambda)$, and therefore both weights vanish. Since the Hilbert--Mumford weights satisfy $\mu^\cL(x,\lambda)=\mu^\cL(\lim_{t \to 0}\lambda(t)x,\lambda)$, it follows that $\mu^\cL(x,\lambda)=0$.
\end{proof}

\begin{proof}[Proof of Theorem~\ref{thm:common-specializations}]
Retain the maximal $F$-split torus $S\subset H$ and the norm $\|\cdot\|$ fixed in Section~\ref{sec:rational-semistability}. We claim that it is enough to prove the theorem in the case when $h_n \in S(F)$. Indeed, by the Cartan decomposition there are compact subsets $C_1,C_2\subset H(F)$ such that $H(F)=C_1S(F)C_2$. We can therefore write $h_n=k_na_n\ell_n$, with $k_n\in C_1$, $a_n\in S(F)$ and $\ell_n\in C_2$. Since $C_1$ and $C_2$ are compact and metrizable, after passing to a subsequence we may assume that $k_n \to k$ and $\ell_n \to \ell$. Thus
\[
\ell_nx_n\longrightarrow\ell x \quad \text{and} \quad a_n\ell_nx_n =k_n^{-1}h_nx_n \longrightarrow k^{-1}x'.
\]

Suppose that the theorem holds for arbitrary sequences of group elements in $S(F)$. Since $\ell x, k^{-1}x' \in X^\ssF$, there exist an $F$-rational one-parameter subgroup $\beta:\Gm\to H$ and an element $b\in H(F)$ such that
\[
\lim_{t\to0}\beta(t)\ell x = \lim_{t\to0}\beta(t)^{-1}bk^{-1}x' \in X^\ssF.
\]
Taking $\lambda=\ell^{-1}\beta\ell$ and $a=\ell^{-1}bk^{-1}$, we obtain
\[
\lim_{t\to0}\lambda(t)x = \ell^{-1}\lim_{t\to0}\beta(t)\ell x = \lim_{t\to0}\lambda(t)^{-1}ax'.
\]
Moreover, $\ell^{-1}\lim_{t\to0}\beta(t)\ell x \in X^\ssF$ by the $H(F)$-invariance of $X^\ssF$. Thus, replacing $x_n$, $x$ and $x'$ by $\ell_nx_n$, $\ell x$ and $k^{-1}x'$, respectively, we may assume from now on that $h_n=a_n\in S(F)$.

Fix a uniformizer $\varpi$ of $F$. Since $S$ is split, we have the split short exact sequence
\[
1 \longrightarrow S(\cO_F) \longrightarrow S(F) \longrightarrow X_*(S) \longrightarrow 0,
\]
and hence
\[
S(F)=S(\cO_F)\,X_*(S)(\varpi),
\]
where $X_*(S)(\varpi) = \{\nu(\varpi): \nu\in X_*(S)\}$. Since $S(\cO_F)$ is compact and metrizable, the same subsequence argument as above allows us to assume that $a_n=\nu_n(\varpi)$ with $\nu_n\in X_*(S)$.

Recall the functions $r_\chi:X^\an\to\bbR$, $\chi \in \wt_S(V)$, defined in \eqref{eq:chi-ratio}. Since the norm $\|\cdot\|$ is split by the $S$-weight decomposition, these ratios transform under $S(F)$ according to
\begin{equation}\label{eq:chi-ratio-torus-action}
r_\chi(ay) = \frac{|\chi(a)|\,r_\chi(y)} {\max_{\eta\in\wt_S(V)}|\eta(a)|\,r_\eta(y)}
\end{equation}
for all $a \in S(F)$ and $y \in X^\an$.

Put $I=\wt_S(x)$ and $J=\wt_S(x')$. For $\chi\in I$ and $\eta\in J$, set $d_n(\chi,\eta) = \langle\chi-\eta,\nu_n\rangle$. Now, take $\chi\in I$ and $\eta\in J$. Since $r_\chi(x_n)\to r_\chi(x)>0$ and $r_\eta(a_nx_n)\to r_\eta(x')>0$, we have $r_\chi(x_n)>0$ and $r_\eta(a_nx_n)>0$ for $n\gg 0$. Formula~\eqref{eq:chi-ratio-torus-action} therefore gives
\[
|\varpi|^{d_n(\chi,\eta)} = \frac{r_\eta(x_n)r_\chi(a_nx_n)} {r_\chi(x_n)r_\eta(a_nx_n)}
\]
for all $n\gg 0$. The denominator tends to the positive number $r_\chi(x) r_\eta(x')$, whereas the numerator tends to $r_\eta(x)r_\chi(x')$, which is positive precisely when $\chi,\eta\in I\cap J$, and is zero otherwise. Since $0<|\varpi|<1$ and $d_n(\chi,\eta)$ takes integer values, we conclude that
\[
\begin{cases}
d_n(\chi,\eta)\text{ is eventually constant}, &\text{if }\chi,\eta\in I\cap J,\\
d_n(\chi,\eta)\longrightarrow+\infty, &\text{otherwise}.
\end{cases}
\]

After discarding finitely many terms, we may assume that $d_n(\chi,\eta)$ is a constant sequence for all $\chi, \eta \in I \cap J$. If $I\neq J$, choose $N$ such that $d_N(\chi,\eta)>d_1(\chi,\eta)$ for all $(\chi,\eta) \in I \times J$ for which $d_n(\chi,\eta)\to +\infty$; if $I=J$, take any $N>1$. Put $\lambda=\nu_N-\nu_1$. By construction, for $(\chi, \eta) \in I \times J$, we have
\begin{equation}\label{eq:chi-eta-separation}
\langle\chi-\eta,\lambda\rangle
\begin{cases}
=0,&\text{if }\chi,\eta\in I\cap J,\\
>0,&\text{otherwise}.
\end{cases}
\end{equation}

In particular, \eqref{eq:chi-eta-separation} implies
\[
\min_{\chi\in I}\langle\chi,\lambda\rangle\ge\max_{\eta\in J}\langle\eta,\lambda\rangle.
\]
On the other hand, formula~\eqref{eq:hilbert-mumford-weight-projective-formula} gives $\mu^\cL(x,\lambda)=-\min_{\chi\in I}\langle\chi,\lambda\rangle$ and $\mu^\cL(x',-\lambda)=\max_{\eta\in J}\langle\eta,\lambda\rangle$. Since $x$ and $x'$ are $S$-semistable, the first of these identities gives $\min_{\chi\in I}\langle\chi,\lambda\rangle\leq0$, while the second gives $\max_{\eta\in J}\langle\eta,\lambda\rangle\ge 0$. Together with the preceding inequality, this forces
\begin{equation}\label{eq:zero-extremal-weights}
\min_{\chi\in I}\langle\chi,\lambda\rangle = \max_{\eta\in J}\langle\eta,\lambda\rangle =0.
\end{equation}

Now, from \eqref{eq:chi-eta-separation} and \eqref{eq:zero-extremal-weights} we obtain
\begin{equation}\label{eq:common-zero-weight-support}
\{\chi\in I:\langle\chi,\lambda\rangle=0\} = I\cap J = \{\eta\in J:\langle\eta,\lambda\rangle=0\}.
\end{equation}
Indeed, by \eqref{eq:zero-extremal-weights} we can fix $\chi_0\in I$ and $\eta_0\in J$ such that $\langle\chi_0,\lambda\rangle=\langle\eta_0,\lambda\rangle=0$. Then $\langle\chi_0-\eta_0,\lambda\rangle=0$, so \eqref{eq:chi-eta-separation} implies that $\chi_0,\eta_0\in I\cap J$. For any $\chi\in I$, applying \eqref{eq:chi-eta-separation} to $(\chi,\eta_0)$ shows that $\langle\chi,\lambda\rangle=0$ if and only if $\chi\in I\cap J$, and similarly for $\eta\in J$ by considering $(\chi_0,\eta)$. This proves~\eqref{eq:common-zero-weight-support} and in particular shows that $I\cap J$ is nonempty.

Put
\[
W=\bigoplus_{\chi\in I\cap J}V_\chi, \qquad W'=\bigoplus_{\chi\notin I\cap J}V_\chi,
\]
and let
\[
\pr_W:\bbP(V)\setminus\bbP(W')\longrightarrow\bbP(W), \quad \left[\sum_{\chi \in \wt_S(V)}v_\chi\right] \longmapsto \left[\sum_{\chi \in I \cap J}v_\chi\right]
\]
be the linear projection. Its Berkovich analytification $\pr_W^\an: \bbP(V)^\an\setminus\bbP(W')^\an\longrightarrow\bbP(W)^\an$ is continuous and $S(F)$-equivariant.

Now, put $K=\cH(x)$ and $K'=\cH(x')$, and represent $x$ and $x'$ by $v\in V_K$ and $v'\in V_{K'}$, respectively. Write $v=\sum_{m\in\bbZ}v_m$ and $v'=\sum_{m\in\bbZ}v'_m$ according to the $\lambda$-weight decomposition of $V$. By \eqref{eq:zero-extremal-weights}, we have $v_m=0$ for $m<0$ with $v_0\neq0$, and $v'_m=0$ for $m>0$ with $v'_0\neq0$. It follows that
\[
[\lambda(t)v]=[v_0 + tv_1 + t^2v_2 + \dots]\longrightarrow [v_0]
\]
in $\bbP(V)(K)$ as $t\to0$ in $F^\times$. Hence $\lim_{t\to0}\lambda(t)x$ is given by the linear projection of $x$ onto $\bbP\left(\bigoplus_{\{\chi\in I:\langle\chi,\lambda\rangle=0\}}V_\chi\right)^\an$. By \eqref{eq:common-zero-weight-support}, the latter is precisely $\bbP(W)^\an$, and therefore
\[
\lim_{t\to0}\lambda(t)x=\pr_W^\an(x).
\]
The same computation in $\bbP(V)(K')$ for $x'$ gives
\[
\lim_{t\to0}\lambda(t)^{-1}x'=\pr_W^\an(x').
\]

It remains to see that the two projections are related by an element of $S(F)$. Since for $\chi,\eta\in I\cap J$ we have $d_n(\chi,\eta)=d_1(\chi,\eta)$, we get $\langle\chi,\nu_n-\nu_1\rangle = \langle\eta,\nu_n-\nu_1\rangle$ for all $n\ge 1$. Consequently, $a_na_1^{-1}$ acts on $W$ by a scalar and therefore acts trivially on $\bbP(W)^\an$. Since the domain where $\pr_W^\an$ is defined is open and contains both $x$ and $x'$, $\pr_W^\an(x_n)$ and $\pr_W^\an(a_nx_n)$ are defined for all sufficiently large $n$. Equivariance of $\pr_W^\an$ and the trivial action of $a_na_1^{-1}$ on $\bbP(W)^\an$ give
\[
\pr_W^\an(a_nx_n)=a_n\pr_W^\an(x_n)=a_1\pr_W^\an(x_n).
\]
Passing to the limit gives $\pr_W^\an(x')=a_1\pr_W^\an(x)$. Therefore, we get
\[
\lim_{t\to0}\lambda(t)x=\pr_W^\an(x)=a_1^{-1}\pr_W^\an(x')=a_1^{-1}\lim_{t\to0}\lambda(t)^{-1}x'= \lim_{t\to0}\lambda(t)^{-1}a_1^{-1}x'.
\]

We conclude by noting that, from \eqref{eq:hilbert-mumford-weight-projective-formula} and \eqref{eq:zero-extremal-weights}, we get $\mu^\cL(x,\lambda)=0$; since $x\in X^\ssF$, Lemma~\ref{lem:semistable-1PS-specialization} therefore implies that $\lim_{t\to 0}\lambda(t)x \in X^\ssF$.
\end{proof}

\section{Cocompactness of weakly admissible period domains}\label{sec:weakly-admissible-period-domains}

In this section we apply Theorem~\ref{thm:semistable-cocompactness} to weakly admissible period domains using the Hilbert--Mumford characterization of weak admissibility of Dat--Orlik--Rapoport \cite[Theorem~9.7.3(ii) and \S9.8]{DOR}. We begin by recalling the definition of the weakly admissible locus and the construction of the corresponding rational linearization.

\subsection{The period-domain setup}

Fix a prime number $p$ and a finite extension $F/\bbQ_p$. Let $\breve F$ be the completion of the maximal unramified extension $F^\mathrm{nr}$ of $F$, and let $\sigma$ denote the Frobenius automorphism of $\breve F/F$. Let $(G,\{\mu\},[b])$ be a period datum over $F$: thus $G$ is a connected reductive group over $F$, $\{\mu\}$ is a $G(\overline F)$-conjugacy class of cocharacters of $G_{\overline F}$, and $[b]\in B(G)$, where $B(G)$ is the Kottwitz set of $\sigma$-conjugacy classes in $G(\breve F)$. We fix a representative $b \in G(\breve F)$ of $[b]$. Let $E=E(G,\{\mu\})$ be the reflex field, and let $\breve E$ be the completion of its maximal unramified extension $E^\mathrm{nr}$. We denote by $\cF(G,\mu)$ the flag variety over $\breve E$ associated to $\{\mu\}$ and by $\cF(G,\mu)^\an$ its Berkovich analytification.

Recall that an \emph{isocrystal} over $F$ is a pair $(D,\varphi)$, where $D$ is a finite-dimensional $\breve F$-vector space and $\varphi:D\to D$ is a $\sigma$-semilinear bijection. For a field extension $K/\breve E$, a \emph{filtered isocrystal} over $K$ is a triple $\cD=(D,\varphi,\Fil^\bullet D_K)$, where $(D,\varphi)$ is an isocrystal over $F$ and $\Fil^\bullet D_K$ is a decreasing, separated, and exhaustive $\bbZ$-filtration on $D_K=D\otimes_{\breve F}K$. Put
\[
t_N(\cD)=\operatorname{val}_F(\det\varphi), \qquad t_H(\cD) = \sum_{i\in\bbZ}i\dim_K\gr^i_{\Fil} D_K.
\]
Here $\operatorname{val}_F$ is normalized by $\operatorname{val}_F(\varpi)=1$ for a uniformizer $\varpi$ of $F$. Since $t_N(\cD)$ does not depend on the filtration, we will also write $t_N(D,\varphi)$, or, suppressing $\varphi$ from the notation, $t_N(D)$. The filtered isocrystal $\cD$ is called weakly admissible if $t_H(\cD)=t_N(\cD)$ and $t_H(\cD')\le t_N(\cD')$ for every filtered subisocrystal $\cD'\subset\cD$ over $K$.

Every $F$-representation $\rho:G\to\GL(W)$ induces a morphism of flag varieties
\[
\rho_\cF:\cF(G,\mu)\longrightarrow\cF(\GL(W),\rho\circ\mu)
\]
over $\breve E$. In particular, for any point $x\in\cF(G,\mu)(K)$, we get a filtration $\Fil_{\rho_\cF(x)}^\bullet$ on $W\otimes_F K$ corresponding to $\rho_\cF(x)\in\cF(\GL(W),\rho\circ\mu)(K)$.

\begin{definition}[{\cite[Definition~1.18]{RZ}}]\label{def:weak-admissibility}
Let $K/\breve E$ be a field extension. A point $x\in\cF(G,\mu)(K)$ is called \emph{weakly admissible} if, for every $F$-representation $\rho:G\to\GL(W)$, the filtered isocrystal
\[
\left(W\otimes_F\breve F,\rho(b)\circ(\id_W\otimes\sigma),\Fil_{\rho_\cF(x)}^\bullet(W\otimes_FK)\right)
\]
is weakly admissible; equivalently, it is enough to test one faithful representation $\rho$.

Let $\cF(G,\mu,b)^{\wa,\mathrm{alg}}\subset\cF(G,\mu)$ denote the set of points $\xi$ for which the canonical $\kappa(\xi)$-valued point is weakly admissible. We then define the weakly admissible locus in $\cF(G,\mu)^\an$ by
\[
\cF(G,\mu,b)^\wa:= \supp_{\cF(G,\mu)}^{-1}\bigl(\cF(G,\mu,b)^{\wa,\mathrm{alg}}\bigr).
\]
\end{definition}

\begin{warning}\label{warn:non-minuscule-weak-admissibility}
For non-minuscule $\mu$, the classical weakly admissible locus considered here must be distinguished from the weakly admissible locus in the $\BdR^+$-affine Grassmannian introduced by Viehmann \cite[Definition~4.1]{Viehmann} for basic $[b]$ and generalized by Chen--Tong \cite[Definition~3.12(2)]{CT25} to arbitrary $[b] \in B(G)$. Indeed, the Bia{\l}ynicki--Birula map is no longer an isomorphism in general (cf.~Section~\ref{sec:admissible-period-domains}), and the two notions of weak admissibility need not agree \cite[Example~4.11]{Viehmann}.
\end{warning}

Let $J_b$ be the $\sigma$-centralizer of $b$, i.e.~the connected reductive $F$-group given by
\[
J_b(R) = \left\{ g\in G(R\otimes_F\breve F) \mathrel{\big|} gb=b\sigma(g) \right\}
\]
for every $F$-algebra $R$, where $\sigma$ acts trivially on $R$. Thus $J_b(F)=\{g\in G(\breve F)\mid gb=b\sigma(g)\}$. The algebraic action of $G_{\breve E}$ on $\cF(G,\mu)$ induces an action of $G(\breve F)$ on $\cF(G,\mu)^\an$, whose restriction to $J_b(F)$ preserves the weakly admissible locus.

The groups $G$ and $J_b$ are defined over $F$, whereas the Hilbert--Mumford characterization of weak admissibility in \cite[Theorem~9.7.3(ii)]{DOR} is formulated over $\bbQ_p$. We reduce to that setting using the restriction-of-scalars construction of \cite[\S9.8]{DOR}. Put
\[
G_0:=\Res_{F/\bbQ_p}G, \qquad H_0:=\Res_{F/\bbQ_p}J_b.
\]
By \cite[Definition~9.8.4 and Example~9.8.5]{DOR}, the element $b$ determines an element $b_0\in G_0(\breve\bbQ_p)$, and there is a canonical isomorphism $H_0\xrightarrow{\sim}J_{b_0}$. In particular, $J_{b_0}(\bbQ_p)=H_0(\bbQ_p)=J_b(F)$.

The conjugacy class $\{\mu\}$ similarly determines a conjugacy class $\{\mu_0\}$ for $G_0$, defined over $E$. The construction of \cite[\S9.8]{DOR} gives an isomorphism of $\breve E$-varieties
\begin{equation}\label{eq:restriction-of-scalars-flag-varieties}
\cF(G,\mu)\overset{\sim}{\longrightarrow}\cF(G_0,\mu_0).
\end{equation}
It follows from \cite[Lemma~9.8.7]{DOR} that under this isomorphism, weak admissibility with respect to $b$ corresponds to weak admissibility with respect to $b_0$.

By \cite[Proposition~9.3.5(i)]{DOR}, together with $H_0\cong J_{b_0}$, there is an isomorphism
\[
H_{0,\breve E}\overset{\sim}{\longrightarrow}Z_{G_{0,\breve E}}(\nu_{b_0}),
\]
where $\nu_{b_0}$ is the Newton cocharacter of $b_0$. Composing this isomorphism with the inclusion of the centralizer into $G_{0,\breve E}$ gives an algebraic action of $H_{0,\breve E}$ on $\cF(G_0,\mu_0)$ and hence, via \eqref{eq:restriction-of-scalars-flag-varieties}, on $\cF(G,\mu)$. After analytification, the induced action of $H_0(\bbQ_p)=J_b(F)$ on $\cF(G,\mu)^\an$ agrees with the action of $J_b(F)$ described above.

\subsection{The linearization for weak admissibility}

We now briefly recall the construction in \cite[\S9.7]{DOR} of the ample rational $H_{0,\breve E}$-linearization on $\cF(G,\mu)$ whose Hilbert--Mumford function detects weak admissibility.

Fix an invariant inner product on $G_0$ in the sense of \cite[\S6.2]{DOR}. Applying the construction of loc.~cit. to $G_0$, we can associate to each conjugacy class $\{\nu\}$ of $\bbQ$-cocharacters of $G_{0,\overline\bbQ_p}$ a $G_0$-equivariant rational line bundle $\cL_\nu$ on the corresponding flag variety $\cF(G_0,\nu)$. More precisely, for a representative $\nu$ over $\overline\bbQ_p$, let $P_\nu\subset G_{0,\overline\bbQ_p}$ be its associated parabolic subgroup. The invariant inner product identifies $\nu$ with a rational character $\nu^*$ of a maximal torus $T\subset G_{0,\overline\bbQ_p}$ through which $\nu$ factors. By invariance of the inner product, $\nu^*$ is trivial on the coroots of the Levi quotient of $P_\nu$ and therefore defines a rational character of $P_\nu$. Over $\overline\bbQ_p$, the associated line bundle is then given by
\[
\cL_\nu= G_{0,\overline\bbQ_p} \times^{P_\nu}\mathbb G_{a,-\nu^*}
\]
on $G_{0,\overline\bbQ_p}/P_\nu$. By \cite[Lemma~6.2.5]{DOR}, it is ample, independent of the choice of representative $\nu\in\{\nu\}$, and descends to the reflex field of $\{\nu\}$ as an equivariant rational line bundle.

Applying this construction to $\{\mu_0\}$ gives an ample $G_{0,\breve E}$-equivariant rational line bundle
\[
\cL_{\mu_0}\in \Pic^{G_{0,\breve E}}(\cF(G_0,\mu_0)) \otimes_\bbZ\bbQ.
\]
Via the identification \eqref{eq:restriction-of-scalars-flag-varieties}, we regard it as a line bundle on $\cF(G,\mu)$.

Put $\lambda_{b_0}:=-\nu_{b_0}$, and let $\cF(G_0,\lambda_{b_0})$ be the associated flag variety over $\breve E$. Similarly, applying the above construction to $\{\lambda_{b_0}\}$ gives an ample $G_{0,\breve E}$-equivariant rational line bundle
\[
\cL_{b_0}\in \Pic^{G_{0,\breve E}}(\cF(G_0,\lambda_{b_0})) \otimes_\bbZ\bbQ.
\]
We regard both $\cL_{\mu_0}$ and $\cL_{b_0}$ as $H_{0,\breve E}$-equivariant rational line bundles by restricting their equivariant structures along $H_{0,\breve E}\xrightarrow{\sim} Z_{G_{0,\breve E}}(\nu_{b_0}) \hookrightarrow G_{0,\breve E}$.

The parabolic subgroup $P_{\lambda_{b_0}}\subset G_{0,\breve E}$ determines a point $y_{b_0}\in\cF(G_0,\lambda_{b_0})(\breve E)$. Since $H_{0,\breve E}$ centralizes $\lambda_{b_0}$, it fixes $y_{b_0}$. Hence the closed immersion
\[
i_{b_0}:\cF(G,\mu)\longrightarrow\cF(G,\mu)\times \cF(G_0,\lambda_{b_0}), \quad x\longmapsto(x,y_{b_0}),
\]
is $H_{0,\breve E}$-equivariant for the diagonal action. Define
\[
\cL_{\mu,b}:=i_{b_0}^*\bigl(\cL_{\mu_0}\boxtimes\cL_{b_0}\bigr) \in \Pic^{H_{0,\breve E}}(\cF(G,\mu))\otimes_\bbZ\bbQ.
\]
Since the pullback of $\cL_{b_0}$ along the constant morphism $x\mapsto y_{b_0}$ is trivial as a rational line bundle, the underlying rational line bundles of $\cL_{\mu,b}$ and $\cL_{\mu_0}$ coincide. In particular, $\cL_{\mu,b}$ is ample.

\subsection{Cocompactness of the weakly admissible locus}

The preceding construction gives the following characterization of weak admissibility.

\begin{proposition}\label{prop:weak-admissibility-criterion}
Let $K/\breve E$ be a field extension and let $x\in\cF(G,\mu)(K)$. Then $x$ is weakly admissible if and only if $\mu^{\cL_{\mu,b}}(x,\lambda)\ge 0$ for every one-parameter subgroup $\lambda:\Gm\to H_0$ defined over $\bbQ_p$. Equivalently, the inequalities may be indexed by the one-parameter subgroups of $J_b$ defined over $F$.
\end{proposition}
\begin{proof}
Since the identification \eqref{eq:restriction-of-scalars-flag-varieties} respects weak admissibility, the first assertion follows by applying \cite[Theorem~9.7.3(ii)]{DOR} to the datum $(G_0,\{\mu_0\},b_0)$ over $\bbQ_p$.\footnote{The standing decency assumption in \cite[\S9.7]{DOR} poses no restriction here, since every class in $B(G_0)$ contains a decent representative \cite[Remark~9.1.34(i)]{DOR}, and the assertion of \cite[Theorem~9.7.3(ii)]{DOR} is invariant under $\sigma$-conjugation of $b_0$.} The second assertion follows from the adjunction $\Hom_{\bbQ_p}\bigl(\mathbb G_{m,\bbQ_p},\Res_{F/\bbQ_p}J_b\bigr) \cong \Hom_F\bigl(\mathbb G_{m,F},J_b\bigr)$.
\end{proof}

Let $\cF(G,\mu)^{\mathrm{ss},\bbQ_p}\subset\cF(G,\mu)^\an$ denote the rational Hilbert--Mumford semistable locus for the action of $H_0$ with linearization $\cL_{\mu,b}$, where semistability is tested against the one-parameter subgroups of $H_0$ defined over $\bbQ_p$.

\begin{corollary}\label{cor:weak-admissibility-semistability}
One has $\cF(G,\mu,b)^\wa = \cF(G,\mu)^{\mathrm{ss},\bbQ_p}$.
\end{corollary}
\begin{proof}
Proposition~\ref{prop:weak-admissibility-criterion} gives $\cF(G,\mu,b)^{\wa,\alg}=\cF(G,\mu)^{\mathrm{ss},\bbQ_p,\alg}$. The assertion now follows by taking inverse images under $\supp_{\cF(G,\mu)}$.
\end{proof}

\begin{theorem}\label{thm:wa-cocompactness}
There exists a compact subset $D\subset\cF(G,\mu,b)^\wa$ such that
\[
J_b(F)\cdot D = \cF(G,\mu,b)^\wa.
\]
\end{theorem}
\begin{proof}
Apply Theorem~\ref{thm:semistable-cocompactness} to the action of the $\bbQ_p$-group $H_0=\Res_{F/\bbQ_p}J_b$ on the projective $\breve E$-variety $\cF(G,\mu)$ with the linearization $\cL_{\mu,b}$. Since $H_0(\bbQ_p)=J_b(F)$, the assertion follows from Corollary~\ref{cor:weak-admissibility-semistability}.
\end{proof}

\begin{remark}[Adic version]\label{rem:adic-wa-cocompactness}
One can similarly consider the weakly admissible locus in the Huber analytification $\cF(G,\mu)^\ad$, defined by
\[
\cF(G,\mu,b)^{\wa,\ad}:=(\supp_{\cF(G,\mu)}^\ad)^{-1}(\cF(G,\mu,b)^{\wa,\alg}).
\]
It follows that $\cF(G,\mu,b)^{\wa,\ad}=\sep_{\cF(G,\mu)}^{-1}(\cF(G,\mu,b)^\wa)$, where $\sep_{\cF(G,\mu)}: \cF(G,\mu)^\ad \to \cF(G,\mu)^\an$ is the adic--Berkovich quotient map. Theorem~\ref{thm:wa-cocompactness} now implies, by the same argument as in Remark~\ref{rem:adic-semistable-cocompactness}, that there exists a quasicompact open subset $U\subset\cF(G,\mu,b)^{\wa,\ad}$ such that $J_b(F)\cdot U=\cF(G,\mu,b)^{\wa,\ad}$.
\end{remark}

\section{Cocompactness and admissible period domains}\label{sec:admissible-period-domains}

In this section we prove Proposition~\ref{prop:admissible-1PS-specialization} about one-parameter specializations in $p$-adic period domains and deduce Theorem~\hyperref[thm:intro-period-domains]{A}(ii). We begin by recalling the description of admissibility in terms of modifications of $G$-bundles on the Fargues--Fontaine curve.

\subsection{Admissibility and modifications of $G$-bundles}

Fix a finite extension $F/\bbQ_p$. For a complete algebraically closed extension $C/F$, let $X_{C^\flat,F}$ denote the Fargues--Fontaine curve over $F$ associated with $C^\flat$. The untilt $C$ determines a point $\infty\in X_{C^\flat,F}$ with residue field $C$ and completed local ring $\widehat{\cO}_{X_{C^\flat,F},\infty}\cong \BdR^+(C)$. When $C$ is understood, we write $\BdR^+$ and $\BdR$ for $\BdR^+(C)$ and $\BdR(C)$, respectively.

The curve $X_{C^\flat,F}$ is complete with respect to the degree function for which every closed point has degree one \cite[Theorem~6.5.2(1)]{FF}. This allows one to define the degree of line bundles on $X_{C^\flat,F}$, which in turn yields the usual notions of degree, slope, and semistability for vector bundles.

For an isocrystal $(D,\varphi)$, let $\cE(D,\varphi)$ be the associated vector bundle on $X_{C^\flat,F}$ constructed in \cite[\S8.2.3]{FF}. We write $\cE(D)$ if $\varphi$ is clear from the context. The functor $(D,\varphi)\mapsto\cE(D,\varphi)$ is an exact tensor functor and sends an isocrystal of slope $\lambda$ to a vector bundle of slope $-\lambda$; see \cite[Proposition~8.2.6]{FF}. In particular, $\deg\cE(D,\varphi)=-t_N(D,\varphi)$.

For a connected reductive group $G/F$, let $\Rep_F(G)$ denote the category of finite-dimensional algebraic $F$-representations of $G$. Recall that an \emph{isocrystal with $G$-structure}, or a \emph{$G$-isocrystal}, is an exact tensor functor from $\Rep_F(G)$ to the category of isocrystals over $F$. Every $b\in G(\breve F)$ defines such a functor by
\[
(W,\rho)\longmapsto \bigl(W\otimes_F\breve F,\, \rho(b)\circ(\id_W\otimes\sigma)\bigr).
\]
Its isomorphism class depends only on the $\sigma$-conjugacy class $[b]$, and this construction induces a bijection between the Kottwitz set $B(G)$ and the set of isomorphism classes of $G$-isocrystals \cite[Remark~3.4(i)]{RR}.

Recall also that a $G$-bundle on $X_{C^\flat,F}$ may equivalently be viewed as an exact tensor functor from $\Rep_F(G)$ to the category of vector bundles on $X_{C^\flat,F}$. Composing the $G$-isocrystal defined by $b$ with the functor $(D,\varphi)\mapsto\cE(D,\varphi)$ gives a $G$-bundle $\cE_b$ on $X_{C^\flat,F}$. By a theorem of Fargues \cite[Th\'eor\`eme~5.1]{Fargues}, the association $[b]\mapsto[\cE_b]$ gives a bijection between $B(G)$ and the set of isomorphism classes of $G$-bundles on $X_{C^\flat,F}$.

The construction of $\cE_b$ gives a canonical trivialization $\tau_b:\cE_b|_{\Spec(\BdR^+)}\xrightarrow{\sim}G_{\Spec(\BdR^+)}$; see \cite[Remark~1.7]{CT26}. Via this trivialization, the $C$-valued points of the $\BdR^+$-affine Grassmannian,
\[
\Gr_G(C)=G(\BdR)/G(\BdR^+),
\]
give rise to modifications of the $G$-bundle $\cE_b$ at $\infty$. Namely, for $z=gG(\BdR^+)$, the corresponding modification $\cE_{b,z}$ is obtained by Beauville--Laszlo gluing \cite{BL} of $\cE_b|_{X_{C^\flat,F}\setminus\{\infty\}}$ with the trivial $G$-bundle over $\Spec(\BdR^+)$, using $g^{-1}\circ\tau_b$ as the gluing isomorphism over $\Spec(\BdR)$.

We now restrict to modifications of a prescribed relative position. Let $\{\mu\}$ be a $G(\overline F)$-conjugacy class of cocharacters of $G_{\overline F}$, and let $E=E(G,\{\mu\})$ be its reflex field. For a complete algebraically closed extension $C/\breve E$, choose an embedding $\overline F\hookrightarrow C$ extending the given embedding $E^{\mathrm{nr}}\hookrightarrow C$. By Hensel's lemma, the embedding $\overline F\hookrightarrow C$ lifts uniquely to $\BdR^+$, thereby giving $\BdR^+$ and its fraction field $\BdR$ an $\overline F$-algebra structure and allowing us to evaluate cocharacters of $G_{\overline F}$ at elements of $\BdR^\times$. The affine Schubert cell associated with $\{\mu\}$ is given on $C$-valued points by
\[
\Gr_{G,\mu}(C) = G(\BdR^+)\mu(\xi)^{-1}G(\BdR^+)/G(\BdR^+) \subset\Gr_G(C),
\]
where $\mu:\mathbb{G}_{m,\overline F}\to G_{\overline F}$ is any representative of $\{\mu\}$ and $\xi$ is a uniformizer of $\BdR^+$. It is independent of the chosen embedding $\overline F\hookrightarrow C$. By \cite[Proposition~3.4.3]{CS}, there is a Bia{\l}ynicki--Birula map
\[
\BB_{G,\mu}(C): \Gr_{G,\mu}(C) \longrightarrow \cF(G,\mu)(C).
\]
To describe this map explicitly, let $P_\mu\subset G_{\overline F}$ be the parabolic subgroup associated with the chosen representative $\mu$, so that $\cF(G,\mu)(C)=G(C)/P_\mu(C)$. Then $\BB_{G,\mu}(C)$ is given by
\begin{equation}\label{eq:BB-map-formula}
g\mu(\xi)^{-1}G(\BdR^+)\longmapsto G(\theta)(g)P_\mu(C),
\end{equation}
where $G(\theta):G(\BdR^+)\to G(C)$ is the group homomorphism induced by the residue map $\theta:\BdR^+\twoheadrightarrow C$.

If $\{\mu\}$ is minuscule, this map is a bijection \cite[Lemma~3.4.4]{CS}; in particular, using $\BB_{G,\mu}(C)^{-1}$ to identify $\cF(G,\mu)(C)$ with $\Gr_{G,\mu}(C)$, every point $x\in\cF(G,\mu)(C)$ determines a modification $\cE_{b,x}$ of $\cE_b$ at $\infty$.

\begin{definition}\label{def:admissibility}
Let $(G,\{\mu\},[b])$ be a period datum over $F$ with $\{\mu\}$ minuscule, and fix a representative $b\in G(\breve F)$ of $[b]$. For $x\in\cF(G,\mu)^\an$, choose a complete algebraically closed extension $C/\cH(x)$ and let $x_C\in\cF(G,\mu)(C)$ be the corresponding $C$-point. We say that $x$ is \emph{admissible} if, for every $F$-representation $\rho:G\to\GL(W)$, the vector bundle $\rho_*\cE_{b,x_C}$ on $X_{C^\flat,F}$ is trivial; equivalently, it is enough to test one faithful representation.
\end{definition}

This definition is independent of the choice of $C$, and we denote by $\cF(G,\mu,b)^\adm$ the resulting locus of admissible points in $\cF(G,\mu)^\an$.

\begin{remark}\label{rem:admissibility-equivalent-definitions}
\begin{enumerate}
\item This Tannakian definition agrees with the definition of admissibility via the triviality of the Harder--Narasimhan point $\nu_{\cE_{b,x_C}}$ used in \cite[Definition~3.1]{CFS}. Indeed, by the functoriality of Harder--Narasimhan points and the fact that a vector bundle on $X_{C^\flat,F}$ is determined up to isomorphism by its Harder--Narasimhan polygon, the bundle $\rho_*\cE_{b,x_C}$ is trivial if and only if $\rho_*(\nu_{\cE_{b,x_C}})=0$. For a faithful representation $\rho$, the latter condition is equivalent to $\nu_{\cE_{b,x_C}}=0$.

\item If $[b]\in B(G,\mu)$, admissibility is also equivalent to the triviality of the $G$-bundle $\cE_{b,x_C}$ \cite[p.~236]{CFS}. If $[b]\notin B(G,\mu)$, however, the bundle $\cE_{b,x_C}$ is never trivial \cite[Proposition~A.9]{Rapoport}, even when $x$ is admissible.
\end{enumerate}
\end{remark}

The locus $\cF(G,\mu,b)^\adm$ is $J_b(F)$-stable and open in $\cF(G,\mu)^\an$, and is contained in the weakly admissible locus $\cF(G,\mu,b)^\wa$. Moreover, the two loci have the same classical points, i.e.~for every finite extension $K/\breve E$, one has $\cF(G,\mu,b)^\adm(K)=\cF(G,\mu,b)^\wa(K)$; see e.g.~\cite[Proposition~3.2]{CFS}.

\subsection{Admissibility under one-parameter specializations}

To study the behavior of admissibility under one-parameter specializations, we would like to relate, for points $x\in\cF(G,\mu)(C)$, the vector bundles $\rho_*\cE_{b,x}$ to the filtered isocrystals corresponding to $\rho_\cF(x)$ that occur in the definition of weak admissibility. The problem is that, even when $\mu$ is minuscule, $\rho\circ\mu$ need not be minuscule, so the corresponding Bia{\l}ynicki--Birula map need not be a bijection, and there is in general no canonical way to assign modifications of vector bundles to points of $\cF(\GL(W),\rho\circ\mu)(C)$. Nevertheless, one can still construct (non-canonical) sections of the Bia{\l}ynicki--Birula map on $C$-valued points that are functorial in the group, which will allow us to recover $\rho_*\cE_{b,x}$ as the modification of $\cE_{\rho(b)}$ along $\rho_\cF(x)$ defined by such a section.

To this end, let $C/\breve E$ be a complete algebraically closed extension. Fix an $\breve F$-algebra section $s:C\to\BdR^+$ of the residue map $\theta:\BdR^+\twoheadrightarrow C$; such a section exists by \cite[Tag~0C34]{Stacks}. Choose an embedding $\overline F\hookrightarrow C$ extending the given embedding $E^{\mathrm{nr}}\hookrightarrow C$, and endow $\BdR^+$ with the $\overline F$-algebra structure given by the composition $\overline F\hookrightarrow C\xrightarrow{s}\BdR^+$.

For any $\{\mu\}$, the choice of $s$ gives a section of $\BB_{G,\mu}(C)$ as follows. For a representative $\mu$ of $\{\mu\}$, let $P_\mu\subset G_{\overline F}$ be the associated parabolic subgroup, so that $\cF(G,\mu)(C)=G(C)/P_\mu(C)$. Applying the functor of points of $G$ to $s$ gives a homomorphism $G(s):G(C)\to G(\BdR^+)$, and for $x=gP_\mu(C)$, we put
\[
\iota_s^{G,\mu}(x):=G(s)(g)\mu(\xi)^{-1}G(\BdR^+),
\]
where $\xi$ is a uniformizer of $\BdR^+$. To see that this is well-defined, we use the defining property of $P_\mu$: for $h\in P_\mu(C)$, the morphism $t\mapsto\mu(t)h\mu(t)^{-1}$ extends to a morphism $f_h:\bbA^1_C\to G_C$. With $\BdR^+$ and $\BdR$ viewed as $C$-algebras via $s$, we have
\[
f_h(\BdR^+)(\xi)=f_h(\BdR)(\xi)=\mu(\xi)G(s)(h)\mu(\xi)^{-1}
\]
in $G(\BdR)$; the first equality follows from functoriality of $f_h$ with respect to the inclusion $\BdR^+\hookrightarrow\BdR$, while the second uses that $\xi\in\BdR^\times$ and that $s$ is an $\overline F$-algebra homomorphism. In particular, $\mu(\xi)G(s)(h)\mu(\xi)^{-1}\in G(\BdR^+)$, so replacing $g$ by $gh$ does not change the coset defining $\iota_s^{G,\mu}(x)$.

The resulting map $\iota_s^{G,\mu}:\cF(G,\mu)(C)\to\Gr_{G,\mu}(C)$ is independent of the choice of the embedding $\overline F\hookrightarrow C$, the representative $\mu$ of $\{\mu\}$, and the uniformizer $\xi$. Moreover, this construction is functorial in $G$. By the description of the Bia{\l}ynicki--Birula map in \eqref{eq:BB-map-formula} and the identity $\theta\circ s=\id_C$, we have
\[
\BB_{G,\mu}(C)\circ\iota_s^{G,\mu}=\id.
\]
In particular, when $\{\mu\}$ is minuscule, $\iota_s^{G,\mu}$ is the inverse of $\BB_{G,\mu}(C)$.

For $\GL(W)$, this section can be described directly in terms of lattices. Putting $D=W\otimes_F\breve F$, the map $\iota_s^{\GL(W),\mu}$ sends a filtration $\Fil^\bullet D_C$ to the $\BdR^+$-lattice
\begin{equation}\label{eq:BB-section-lattice}
\Lambda_s(D,\Fil^\bullet D_C) = \sum_{m\in\bbZ}\xi^{-m}\bigl(\Fil^mD_C\otimes_{C,s}\BdR^+\bigr)\subset D\otimes_{\breve F}\BdR.
\end{equation}
For any filtered isocrystal $(D,\Fil^\bullet D_C)$, we denote by $\cE_s(D,\Fil^\bullet D_C)$ the modification of $\cE(D)$ obtained by Beauville--Laszlo gluing with this lattice.

We will use the compatibility of these modifications with strict exact sequences of filtered isocrystals. A sequence
\[
0\longrightarrow(D',\Fil')\longrightarrow(D,\Fil)\overset{q}{\longrightarrow}(D'',\Fil'')\longrightarrow 0
\]
is \emph{strict exact} if it is exact as a sequence of isocrystals and, for every $m\in\bbZ$, the sequence
\[
0\longrightarrow\Fil'^mD'_C\longrightarrow\Fil^mD_C\overset{q_C}{\longrightarrow}\Fil''^mD''_C\longrightarrow 0
\]
is exact. For such a sequence, strictness implies that the induced sequence of filtered $C$-vector spaces splits. Applying formula \eqref{eq:BB-section-lattice} to such a splitting, and using that $D_C\otimes_{C,s}\BdR^+\cong D\otimes_{\breve F}\BdR^+$ naturally in the isocrystal $D$, we obtain
\[
\Lambda_s(D',\Fil')=\Lambda_s(D,\Fil)\cap(D'\otimes_{\breve F}\BdR),\qquad \Lambda_s(D'',\Fil'')=(q\otimes 1)(\Lambda_s(D,\Fil)).
\]
By \cite[Remark~2.3]{NV}, these identities imply that the induced sequence
\[
0\longrightarrow\cE_s(D',\Fil') \longrightarrow\cE_s(D,\Fil) \longrightarrow\cE_s(D'',\Fil'') \longrightarrow0
\]
of vector bundles on $X_{C^\flat,F}$ is exact. In particular, taking $\cE_s(\cdot)$ commutes with finite direct sums of filtered isocrystals. For any filtered isocrystal $(D,\Fil)$, we shall also use the degree formula
\begin{equation}\label{eq:modification-degree}
\deg\cE_s(D,\Fil) = t_H(D,\Fil)-t_N(D);
\end{equation}
see \cite[\S2.2, p.~231]{CFS}.

For the remainder of this section, fix a period datum $(G,\{\mu\},[b])$ over $F$ and a representative $b\in G(\breve F)$ of $[b]$.

\begin{proposition}\label{prop:admissible-1PS-specialization}
Assume that $\{\mu\}$ is minuscule. Let $x\in\cF(G,\mu,b)^\adm$ and let $\lambda:\Gm\to J_b$ be a one-parameter subgroup over $F$. Let $x_0=\lim_{t\to0}\lambda(t)x\in\cF(G,\mu)^\an$. If $x_0\in\cF(G,\mu,b)^\wa$, then $x_0\in\cF(G,\mu,b)^\adm$.
\end{proposition}
\begin{proof}
Since the limit defining $x_0$ can be computed in $\cF(G,\mu)(\cH(x))$, we have a canonical embedding $\cH(x_0)\hookrightarrow\cH(x)$. Choose a complete algebraically closed extension $C/\cH(x)$, and denote again by $x,x_0$ the corresponding $C$-valued points of $\cF(G,\mu)$. With this notation, $\lim_{t\to0}\lambda(t)x$ is now understood as a limit in $\cF(G,\mu)(C)$.

Let $\rho:G\to\GL(W)$ be an $F$-representation, and put $(D,\varphi):=(W\otimes_F\breve F,\rho(b)\circ(\id_W\otimes\sigma))$. The representation $\rho$ induces a homomorphism $J_b\to J_{\rho(b)}$, where $J_{\rho(b)}$ is the algebraic automorphism group of the isocrystal $(D,\varphi)$. Let $\rho(\lambda):\Gm\to J_{\rho(b)}$ be the composition of $\lambda$ with this homomorphism.

Fix an $\breve F$-algebra section $s:C\to\BdR^+$ of the residue map as above. By the functoriality of the section $\iota_s^{G,\mu}$, the representation $\rho$ gives a commutative diagram
\[
\begin{tikzcd}[column sep=5em,row sep=3em]
\cF(G,\mu)(C) \arrow[r,"\iota_s^{G,\mu}","\cong"'] \arrow[d,"\rho_\cF(C)"']
& \Gr_{G,\mu}(C) \arrow[d,"\rho_{\Gr}(C)"] \\
\cF(\GL(W),\rho\circ\mu)(C) \arrow[r,"\iota_s^{\GL(W),\rho\circ\mu}"']
& \Gr_{\GL(W),\rho\circ\mu}(C)
\end{tikzcd}
\]
where $\rho_{\Gr}$ is the map induced by $\rho$ on the $\BdR^+$-affine Grassmannians. The commutativity of this diagram and the functoriality of Beauville--Laszlo gluing give
\[
\cE_s(D,\Fil_{\rho_\cF(y)})\cong\rho_*\cE_{b,y}
\]
for every $y\in\cF(G,\mu)(C)$. Since $x$ is admissible, the vector bundle $\cE_s(D,\Fil_{\rho_\cF(x)})$ is thus trivial. We need to show that $\cE_s(D,\Fil_{\rho_\cF(x_0)})$ is also trivial; since $\rho$ is arbitrary, this will prove the admissibility of $x_0$.

By equivariance and continuity of $\rho_\cF(C)$, we have
\[
\rho_\cF(x_0)=\lim_{t\to0}\rho(\lambda)(t)\cdot\rho_\cF(x)
\]
in $\cF(\GL(W),\rho\circ\mu)(C)$. Let $D=\bigoplus_{i\in\bbZ}D_i$ be the weight space decomposition according to the action of $\rho(\lambda)$, and set $U_i=\bigoplus_{m\ge i}D_m$. Since $\rho(\lambda)$ acts by automorphisms of the isocrystal $(D,\varphi)$, the subspaces $D_i$ and $U_i$ are subisocrystals of $D$. The calculation in the proof of \cite[Lemma~2.2.2]{DOR} shows that the limit $\lim_{t\to0}\rho(\lambda)(t)\cdot\rho_\cF(x)$ in $\cF(\GL(W),\rho\circ\mu)(C)$ is given exactly by the associated graded filtration $\gr_U\Fil_{\rho_\cF(x)}$ on $D_C$. Namely, for $m \in \bbZ$, the $m$-th piece of $\gr_U\Fil_{\rho_\cF(x)}$ is given by
\[
(\gr_U\Fil_{\rho_\cF(x)})^m = \bigoplus_{i \in \bbZ} \frac{\Fil_{\rho_\cF(x)}^m D_C \cap U_{i,C}}{\Fil_{\rho_\cF(x)}^m D_C \cap U_{i+1,C}},
\]
where the $i$-th summand is naturally a subspace of $D_{i,C}$. Therefore,
\begin{equation}\label{eq:associated-graded-filtration-splitting}
(D,\Fil_{\rho_\cF(x_0)}) = \bigoplus_{i \in \bbZ}(D_i,\Fil_{D_i})
\end{equation}
as filtered isocrystals, where $\Fil_{D_i}$ is the filtration on $D_{i,C}$ whose $m$-th piece is given by
\[
\Fil_{D_i}^m = \frac{\Fil_{\rho_\cF(x)}^m D_C \cap U_{i,C}}{\Fil_{\rho_\cF(x)}^m D_C \cap U_{i+1,C}}.
\]

For each $i\in\bbZ$, equip $U_i$ with the filtration $\Fil_{U_i}$ defined by $\Fil_{U_i}^m=\Fil_{\rho_\cF(x)}^m D_C\cap U_{i,C}$. Then the sequence of filtered isocrystals
\[
0\longrightarrow (U_{i+1},\Fil_{U_{i+1}}) \longrightarrow (U_i,\Fil_{U_i}) \longrightarrow (D_i,\Fil_{D_i}) \longrightarrow 0
\]
is strict exact. By the compatibility established above, applying $\cE_s(\cdot)$ gives a short exact sequence
\begin{equation}\label{eq:bundle-SES}
0\longrightarrow\cV_{i+1} \longrightarrow \cV_i \longrightarrow \cQ_i \longrightarrow 0
\end{equation}
of vector bundles on $X_{C^\flat,F}$, where $\cV_i=\cE_s(U_i,\Fil_{U_i})$ and $\cQ_i=\cE_s(D_i,\Fil_{D_i})$. Note that $\cV_i=\cE_s(D,\Fil_{\rho_\cF(x)})$ for $i\ll 0$, while $\cV_i=0$ for $i\gg 0$.

Applying $t_H-t_N$ to \eqref{eq:associated-graded-filtration-splitting} gives
\[
t_H(D,\Fil_{\rho_\cF(x_0)})-t_N(D)= \sum_i\bigl(t_H(D_i,\Fil_{D_i})-t_N(D_i)\bigr).
\]
The filtered isocrystal $(D,\Fil_{\rho_\cF(x_0)})$ is weakly admissible, so the left-hand side vanishes, while every summand on the right is $\leq0$. Consequently, $t_H(D_i,\Fil_{D_i})-t_N(D_i)=0$ for all $i$, and \eqref{eq:modification-degree} gives $\deg\cQ_i=0$ for all $i$. Taking degrees in \eqref{eq:bundle-SES}, starting from $\cV_i=0$ for $i\gg0$, we obtain $\deg\cV_i=0$ for every $i\in\bbZ$.

Each $\cV_i$ is a degree-zero subbundle of the trivial vector bundle $\cE_s(D,\Fil_{\rho_\cF(x)})$, and is therefore trivial. It follows from \eqref{eq:bundle-SES} that all the $\cQ_i$ are then trivial as well. Finally, \eqref{eq:associated-graded-filtration-splitting} and compatibility of $\cE_s(\cdot)$ with direct sums give $\cE_s(D,\Fil_{\rho_\cF(x_0)})\cong\bigoplus_i\cQ_i$, which is therefore trivial.
\end{proof}

\subsection{Admissibility and orbit closures}

Theorem~\ref{thm:common-specializations} and Proposition~\ref{prop:admissible-1PS-specialization} now imply the following:

\begin{corollary}\label{cor:admissible-orbit-limits}
Assume that $\mu$ is minuscule. Suppose that
\[
x_n\longrightarrow x,\qquad h_nx_n\longrightarrow x'
\]
in $\cF(G,\mu)^\an$, where $x_n\in\cF(G,\mu)^\an$, $h_n\in J_b(F)$, $x\in\cF(G,\mu,b)^\adm$ and $x'\in\cF(G,\mu,b)^\wa$. Then $x'\in\cF(G,\mu,b)^\adm$.
\end{corollary}
\begin{proof}
Since admissibility implies weak admissibility, both $x$ and $x'$ are weakly admissible, and hence semistable for the $\bbQ_p$-group $H_0$ by Corollary~\ref{cor:weak-admissibility-semistability}. Applying Theorem~\ref{thm:common-specializations} to this action, and using the equality $H_0(\bbQ_p)=J_b(F)$ and the cocharacter correspondence in Proposition~\ref{prop:weak-admissibility-criterion}, we obtain an $F$-rational cocharacter $\lambda:\Gm\to J_b$ and $a\in J_b(F)$ such that
\[
y:=\lim_{t\to0}\lambda(t)x =\lim_{t\to0}\lambda(t)^{-1}ax'\in\cF(G,\mu,b)^\wa.
\]
By Proposition~\ref{prop:admissible-1PS-specialization}, $y\in\cF(G,\mu,b)^\adm$. Let $\varpi$ be a uniformizer of $F$. Since $\lambda(\varpi^n)^{-1}ax'\to y$ and $\cF(G,\mu,b)^\adm$ is open, we have $\lambda(\varpi^n)^{-1}ax'\in\cF(G,\mu,b)^\adm$ for all sufficiently large $n$. The $J_b(F)$-invariance of $\cF(G,\mu,b)^\adm$ then gives $x'\in\cF(G,\mu,b)^\adm$.
\end{proof}

Since $\breve E$ has a countable dense subset, the space $\cF(G,\mu)^\an$ is metrizable \cite[Remark~1.5]{HLP}. In particular, closedness in $\cF(G,\mu)^\an$ and its subspaces is equivalent to sequential closedness, while compactness is equivalent to sequential compactness. We obtain the following consequence for compact subsets of the admissible locus.

\begin{theorem}\label{thm:admissible-orbit-closure}
Assume that $\mu$ is minuscule. For every compact subset $C\subset\cF(G,\mu,b)^\adm$, one has
\[
\overline{J_b(F)\cdot C}\cap\cF(G,\mu,b)^\wa\subset\cF(G,\mu,b)^\adm.
\]
\end{theorem}
\begin{proof}
Let $x' \in \overline{J_b(F)\cdot C}\cap\cF(G,\mu,b)^\wa$. By metrizability, there are $x_n\in C$ and $h_n\in J_b(F)$ such that $h_nx_n\to x'$. After passing to a subsequence, compactness of $C$ gives $x_n\to x\in C \subset \cF(G,\mu,b)^\adm$. Corollary~\ref{cor:admissible-orbit-limits} now implies $x'\in\cF(G,\mu,b)^\adm$.
\end{proof}

We can now prove the conjecture of Chen--Fargues--Shen \cite[Conjecture~7.2]{CFS}. The conjecture was originally stated only for basic $[b] \in B(G,\mu)$. Our proof shows that it in fact holds for arbitrary $[b] \in B(G)$, without the basic assumption.

\begin{theorem}\label{thm:cfs-equivalence}
Assume that $\mu$ is minuscule. Then, for every $[b]\in B(G)$, the following conditions are equivalent:
\begin{enumerate}
\item $\cF(G,\mu,b)^\adm=\cF(G,\mu,b)^\wa$;
\item there is a compact subset $D\subset\cF(G,\mu,b)^\adm$ with $J_b(F)\cdot D=\cF(G,\mu,b)^\adm$.
\end{enumerate}
\end{theorem}
\begin{proof}
The implication (i)$\Rightarrow$(ii) follows from Theorem~\ref{thm:wa-cocompactness}. Conversely, let $D$ be as in (ii). Since the loci $\cF(G,\mu,b)^\adm$ and $\cF(G,\mu,b)^\wa$ have the same classical points, which are dense in $\cF(G,\mu,b)^\wa$ \cite[\S3.3, p.~57]{Berkovich}, we get $\overline{J_b(F)\cdot D}\cap\cF(G,\mu,b)^\wa = \overline{\cF(G,\mu,b)^\adm}\cap\cF(G,\mu,b)^\wa = \cF(G,\mu,b)^\wa$. Theorem~\ref{thm:admissible-orbit-closure} therefore gives $\cF(G,\mu,b)^\wa = \cF(G,\mu,b)^\adm$.
\end{proof}

\begin{remark}[Adic version]\label{rem:adic-cfs-equivalence}
Theorem~\ref{thm:cfs-equivalence} also holds in the adic setting. Indeed, assume that $\mu$ is minuscule and let $\cF(G,\mu,b)^{\adm,\ad}$ and $\cF(G,\mu,b)^{\wa,\ad}$ denote the admissible and weakly admissible loci in the adic analytification $\cF(G,\mu)^\ad$. Since they satisfy $\cF(G,\mu,b)^{\adm,\ad}=\sep_{\cF(G,\mu)}^{-1}(\cF(G,\mu,b)^\adm)$ and $\cF(G,\mu,b)^{\wa,\ad}=\sep_{\cF(G,\mu)}^{-1}(\cF(G,\mu,b)^\wa)$, we get
\[
\cF(G,\mu,b)^\adm=\cF(G,\mu,b)^\wa \quad\Longleftrightarrow\quad \cF(G,\mu,b)^{\adm,\ad}=\cF(G,\mu,b)^{\wa,\ad}.
\]
Similarly, condition (ii) of Theorem~\ref{thm:cfs-equivalence} is equivalent to the existence of a quasicompact open subset $U\subset \cF(G,\mu,b)^{\adm,\ad}$ such that $J_b(F)\cdot U=\cF(G,\mu,b)^{\adm,\ad}$. Indeed, given such a subset $U$, take $D=\sep_{\cF(G,\mu)}(U)$. Conversely, given $D \subset \cF(G,\mu,b)^\adm$ as in (ii), take $U$ to be any quasicompact open subset of $\cF(G,\mu,b)^{\adm,\ad}$ containing $\sep_{\cF(G,\mu)}^{-1}(D)$.
\end{remark}

\bibliographystyle{amsalpha}
\bibliography{bibliography}
\end{document}